\documentclass[12pt]{amsart}

\newsymbol\pp 1275
\newsymbol\twoheadrightarrow 1310

\newtheorem{theorem}{Theorem}
\newtheorem{proposition}[theorem]{Proposition}

\newtheorem{lemma}[theorem]{Lemma}
\newtheorem{corollary}[theorem]{Corollary}
\newcommand{\wt}{\operatorname{wt}}

\theoremstyle{definition}

\newtheorem{definition}[theorem]{Definition}

\newcommand{\B}{{B}}

\newcommand{\Z}{{\mathbb Z}}

\newcommand{\F}{{\mathbb F}}

\title{Error Correcting Codes and $B_h$-Sequences}
\author{Harm Derksen}
\thanks{The author was supported by NSF, grant DMS 0102193.}

\begin{document}
\maketitle
\begin{abstract}
We construct 
error-correcting (nonlinear) binary codes using a construction
of Bose and Chowla in additive number theory. Our method
 extends a
construction of Graham and Sloane for constant
weight codes. The new codes improve 1028 of the 7168
best known $h$-error correcting codes of wordlength $\leq 512$ with $1\leq h\leq 14$.
We give asymptotical comparisons to shortened BCH codes.
\end{abstract}
\noindent {\it Keywords: $A(n,d)$,
$A(n,d,w)$, binary codes, $B_h$-sequences, error correcting, tables.}
\section{Introduction}
\noindent Let $\F_2$ be the field with two elements.
The following standard definitions are of particular interest to
coding theory:
\begin{eqnarray*}
A(n,d) &:= & \max\{\#C\mid \mbox{$C\subseteq \F_2^n$
is a code with minimum distance $\geq d$}\},\\
A(n,d,w) &:=& \max\{\#C\mid \mbox{$C\subseteq \F_2^n$ has minimum distance $\geq d$
and constant weight $w$}\}.
\end{eqnarray*}
The goal of this paper is to use $B_h$-sequences to construct good codes.
For various values of $n$ and $d$ we will improve the best known
lower bound for $A(n,d)$.

A table of best known lower bounds for $A(n,d)$ for $n\leq 512$ and $d\leq 29$
can be found on the webpages maintained by Litsyn, Rains, and 
Sloane (see~\cite{LRS}). These webpages are based on the published
tables of best known codes
by MacWilliams and Sloane (see~\cite{McWS}) and
more recently by  Litsyn (see~\cite{Li}).
A table with lower bounds for $A(n,d,w)$ can be found in~\cite{McWS}
or at the webpage~\cite{RS}.
For many values of $n$ and $d$ the constructions in this paper give
better lower bounds for $A(n,d)$ than known up to now.
We supplied a table of all improved lower bounds for $A(n,d)$ in
Section~\ref{sec6}.

In Section~\ref{sec2} we discuss
$\B_h$-sequences and their relation to error-correcting codes.
In Section~\ref{sec3} we present a  generalized
Bose-Chowla construction for $B_h$-sequences.
 We apply
the construction of $B_h$-sequences to 
obtain good lower bounds for $A(n,d)$.
A similar construction was already used by Graham and Sloane
 for constant weight codes (see~\cite{GS}). It seems that up to now
it has been overlooked that this construction also improves 
lower bounds for non-constant weight codes.
In Section~\ref{sec4} we study the asymptotics 
of our bounds and we compare our bounds to the bounds
derived from shortened BCH codes.

\noindent{\bf Acknowledgement.}
I would like to thank Sergey Fomin, Simon Litsyn and Joachim Rosenthal
for helpful discussions and suggestions, and Neil Sloane for
pointing out the paper \cite{GS} to me. 
\section{$\B_h$-sequences}\label{sec2}
\noindent The notion of $\B_h$-sequence comes from additive number theory
(see for example \cite[Chapter II]{HR}).
\begin{definition}
A sequence $g_1,\dots,g_n$ in an abelian group $G$ is called
a $B_h$-sequence if all
$$
g_{i_1}+g_{i_2}+\cdots+g_{i_h},\quad 
1\leq i_1\leq i_2\leq\cdots\leq i_h\leq n
$$ 
are distinct.
\end{definition}
\noindent We now discuss the relation between $B_h$-sequences
and codes. For two words $x,y\in \F_2$ we will write
$\wt(x)$ for the weight of $x$ and $\delta(x,y)=\wt(x-y)$ for the
Hamming distance between $x$ and $y$.
\begin{proposition}\label{prop1}
 Suppose that $g_1,\dots,g_n$ is a $B_h$-sequence in
an abelian group $G$.
Identify $\F_2$ with $\{0,1\}\subset \Z$ and define the map
$\phi:\F_2^n\to G$ by
$$
\phi(x_1,\dots,x_n)=x_1g_1+x_2g_2+\cdots+x_ng_n.
$$ 
 For every integer $w$ and every $g\in G$ we have
that 
$$
\phi^{-1}(g)\cap \{x\in \F_2^n\mid \wt(x)=w\}
$$
is a code with minimum distance $\geq 2h+2$ for which all codewords
have weight $w$.
\end{proposition}
\begin{proof}
 Suppose that $x,y\in \F_2^n$ such that $\wt(x)=\wt(y)=w$,
$\phi(x)=\phi(y)=g\in G$ and
$\delta(x,y)<2h+2$. We will show that $x=y$.
In $G$ we have the equality
\begin{equation}\label{eq10}
x_1g_1+x_2g_2+\cdots+x_ng_n=y_1g_1+y_2g_2+\cdots+y_ng_n
\end{equation}
Let $X=\{i\mid x_i=1, y_i=0\}$ and $Y=\{i\mid x_i=0,y_i=1\}$.
Note that $\#X=\#Y$ because $\wt(x)=\wt(y)$ and
$\#X+\#Y=\delta(x,y)<2h+2$. In particular $\delta(x,y)$ is
even and $\#X=\#Y\leq h$.
Now (\ref{eq10}) simplifies to
$$
\sum_{i\in X}g_i=\sum_{i\in Y}g_i.
$$
We have
$$
(h-\#X)g_1+\sum_{i\in X}g_i=(h-\#Y)g_1+\sum_{i\in Y}g_i.
$$
We conclude that
 $X=Y$ because $g_1,g_2,\dots,g_n$ is a $B_h$-sequence.
Since $X$ and $Y$
are obviously disjoint we get $X=Y=\emptyset$
and therefore $x=y$.
\end{proof}
\begin{definition}
Let $c(n,h)$ be the smallest positive integer $C$ such
that an abelian group $G$ with cardinality $C$ and a $B_h$-sequence
of cardinality $n$ in $G$ exist.
\end{definition}
\begin{corollary}\label{cor1}
We have the following inequality:
$$
A(n,2h+2,w)\geq \frac{1}{c(n,h)}{n\choose w}.
$$
\end{corollary}
\begin{proof}
Suppose that $g_1,g_2,\dots,g_n$ is a $B_h$-sequence
in an abelian group $G$ with $\#G=c(n,h)$.
Let $\phi:\F_2^n\to G$ be as in Proposition~\ref{prop1}. 
 By the pigeonhole principle,
at least for one $g\in G$ we must have
$$
\#\{x\in \F_2^n\mid \wt(x)=w\}\cap \phi^{-1}(g)\geq 
\frac{1}{\#G}\#\{x\in \F_2^n\mid \wt(x)=w\}=\frac{1}{\#G}{n\choose w}
$$
and by Proposition~\ref{prop1} we have
that
$$
\{x\in \F_2^n\mid \wt(x)=w\}\cap \phi^{-1}(g)
$$
is a code consisting of words with weight $w$ having minimum distance $2h+2$.
\end{proof}

\section{A construction for $B_h$ sequences}\label{sec3}
Let $q$ be a prime power and let 
$\F_q$ be the field with $q$ elements.
Bose and Chowla constructed $B_h$-sequences in the cyclic group $\F_q^\star$
and the integers (see~\cite{BC} or
{\cite[II~\S 2, Theorem 3]{HR}}).
Suppose that $P(X)$ is a polynomial with coefficients in $\F_q$.
We will work in the ring $\F_q[X]/P(X)$. The multiplicative
group of invertible elements in  $\F_q[X]/P(X)$ is denoted by $(\F_q[X]/P(X))^\star$.
The construction of Bose and Chowla easily generalizes to
the construction of $B_h$-sequences in the (not
necessarily cyclic) group $(\F_q[X]/P(X))^\star$.
Note that $\F_q^\star$, the multiplicative group of units in $\F_q$, is a
subgroup of $(\F_q[X]/P(X))^\star$. 
\begin{theorem}\label{theo1}
Suppose that $P(X)$ is a polynomial in $\F_q[X]$, and
$a_1,\dots,a_n\in \F_q$ are distinct  with
$P(a_i)\neq 0$ for all $i$.
\begin{enumerate}
\renewcommand{\theenumi}{\alph{enumi}}
\item If  $P(X)$ has degree $h$, then
the sequence 
$$X-a_1,X-a_2,\cdots,X-a_n$$ 
is a $\B_h$-sequence
in the (multiplicative) group $(\F_q[X]/P(X))^\star$.
\item Suppose that $P(X)$ has degree $h+1$ then
(the image of) 
$$1,X-a_1,X-a_2,\cdots,X-a_n$$ 
is a $\B_h$-sequence in $(\F_q[X]/P(X))^\star/\F_q^\star$.
\end{enumerate}
\end{theorem}
\begin{proof}
The proof is similar to the original proof of Bose and 
Chowla (see~\cite{BC} or
{\cite[II~\S 2, Theorem 3]{HR}}).
\end{proof}
We would like to give upper bounds for $c(n,h)$
using the previous theorem.
Let us define $\mu(P(X))=\#(\F_q[X]/P(X))^\star$
for any polynomial $P(X)\in \F_q[X]$.
Let us first see how we can minimalize $\mu(P(X))$
by varying the choice of $P(X)$.
\begin{lemma}\label{lem51}
Let $S\subseteq \F_q$ be a subset with $\#S=n$
and suppose that $h$ is a positive integer
such that $h<q$. The smallest
possible value of $\mu(P(X))$
where $P(X)$ is a polynomial of degree $h$ satisfying
$P(x)\neq 0$ for all $x\in S$ will be denoted by
$\mu(q,n,h)$. The values of $\mu(q,n,h)$ are given
in the following table:\\

\hfill\begin{tabular}{||l|l||}
\hline\hline
$q,n,h$ & $\mu(q,n,h)$\\
\hline
$q\geq n+h$ & $(q-1)^h$ \\
$n<q< n+h$, $n+h-q$ even & $(q-1)^{q-n}(q^2-1)^{(n+h-q)/2}$\\
$n<q<n+h$, $n+h-q$ odd & $q(q-1)^{q-n}(q^2-1)^{(n+h-q-1)/2}$\\
$q=n$, $h$ even & $(q^2-1)^{h/2}$\\
$q=n$, $h$ odd & $(q^2-1)^{(h-3)/2}(q^3-1)$\\
\hline\hline
\end{tabular}\hfill\break
\end{lemma}
\begin{proof}
Suppose that $r\geq 1$ is a positive integer.
Note that $\F_{q^r}$ has at least $q$ elements that
do not lie in any proper subfield. For example, one can
take $\alpha,\alpha+1,\dots,\alpha+q-1$ where
$\alpha\in \F_{q^r}$ is a generator over the field
$\F_q$. This implies that if $P(X)$ is a polynomial of degree $<q$
and $r$ is a positive integer, then one can find an irreducible
polynomial $Q(X)\in \F[X]$ of degree $r$ such that $P(X)$ and $Q(X)$
are relatively prime. We will use this repeatedly in the proof.

If $P(X),Q(X)\in \F_q[X]$ are polynomials and $Q(X)$ has degree $d$, then
$$
\mu(P(X))\mu(Q(X))\leq \mu(P(X)Q(X))\leq \mu(P(X))q^d.
$$
The left inequality  is an equality
 if and only if $P(X)$ and $Q(X)$ are relatively prime. The right
 inequality is an equality if and only if every irreducible  factor
 of $Q(X)$ is an irreducible factor of $P(X)$.

Suppose now that $P(X)$ is a polynomial as in the Lemma
such that $\mu(P(X))$ is minimal.
First assume that $S\neq \F_q$. 

Suppose $P(X)$ has
an irreducible factor $P_1(X)$ of degree $d\geq 3$.
We can write $P(X)=P_1(X)P_2(X)$.
Let $b\in \F_q\setminus S$ and let $P_3(X)$ be an
irreducible polynomial of degree $d-1$ relatively
prime to $P_2(X)$.

 Then we get
$$
\mu(P(X))\geq \mu(P_1(X))\mu(P_2(X))=
(q^{d}-1)\mu(P_2(X))>
(q^{d-1}-1)q\mu(P_2(X)) \geq
$$
$$
\geq\mu(P_3(X))q\mu(P_2(X))\geq
\mu(P_3(X))\mu((X-b)P_2(X))=\mu(P_3(X)(X-b)P_2(X)).
$$
This shows that $P(X)$ can only have linear and quadratic factors.

Suppose that $P(X)$ has a irreducible quadratic factor $P_1(X)$
of higher multiplicity. We can write $P(X)=P_1(X)^2P_2(X)$.
Let $P_3(X)$ be an irreducible quadratic polynomial which
is relatively prime to $P(X)$. Then we have
$$
\mu(P(X))=q^2\mu(P_1(X)P_2(X))>(q^2-1)\mu(P_1(X)P_2(X))=\mu(P_1(X)P_2(X)P_3(X)).
$$
This is a contradiction, so all quadratic factors of $P(X)$ should
be distinct. 

If $P(X)$ has a quadratic factor then
$P(X)$ must vanish on all $x\in \F_q\setminus S$ because
otherwise we could replace the quadratic factor by $(X-b)^2$
where $X-b$ is relative prime to $P(X)$ and get a smaller
value of $\mu(P(X))$. In a similar way,
it also follows that if $P(X)$ does not vanish on all $x\in \F_q\setminus S$,
then $P(X)$ cannot have multiple zeroes in $\F_q$ because
otherwise we could replace a factor $(X-b)^2$ of $P(X)$ 
by $(X-b)(X-c)$ where $X-c$ is relatively prime to $P(X)$.
If
$P(X)$ has a zero of multiplicity $\geq 3$,
say $(X-b)^3$ divides $P(X)$, then
we could replace the factor $(X-b)^2$
by an irreducible quadratic polynomial of
degree 2 and decrease $\mu(P(X))$.
This shows that $P(X)$ has zeroes in $\F_q$ with
multiplicities at most 2. Suppose that $P(X)$
has two distinct zeroes in $\F_q$ with multiplicity 2,
say $(X-b)^2(X-c)^2$ divides $P(X)$. We
can replace the factor $(X-b)(X-c)$ by an irreducible
quadratic polynomial relatively prime to $P(X)$ which
would decrease $\mu(P(X))$ again.
This shows that $P(X)$ has at most one zero in $\F_q$
with multiplicity 2. 

If $q-n\geq h$ then $P(X)$ cannot have quadratic factors
or multiple zeroes, because $P(X)$
cannot vanish on $\F_q\setminus S$. In that case $P(X)$ is a product
of distinct linear factors and $\mu(P(X))=(q-1)^h$.
If $q-n<h$, then $P(X)$ must have a quadratic factor or a multiple zero
in $\F_q$. Therefore $P(X)$ must vanish on $\F_q\setminus S$.
Since there is at most one zero with multiplicity 2, there
are either $q-n$ or $q-n+1$ linear factors.
If $h+n-q$ is even, then there are exactly $q-n$
distinct linear factors and $(h+n-q)/2$ distinct quadratic
irreducible factors. In that case we have
$\mu(P(X))=(q-1)^{q-n}(q^2-1)^{(h+n-q)/2}$.
If $h+n-q$ is odd, then $P(X)$ has $q-n$ zeroes
of which one has multiplicity 2. The number
of quadratic irreducible factors is $(h+n-q-1)/2$.
In this case we have
$\mu(P(X))=q(q-1)^{q-n}(q^2-1)^{(h+n-q-1)/2}$.

Now we need to discuss the case that $S=\F_q$. Similar arguments
as before show that now $P(X)$ cannot have factors of 
degree $\geq 4$ and that all factors of $P(X)$
are distinct.
Now $P(X)$ can have at most 1 irreducible factor of degree 3
because otherwise two irreducible factors of degree 3 could
be replaced by 3 distinct irreducible factors of degree 2 
which are relatively prime to $P(X)$.
If $h$ is even, then $P(X)$ has $h/2$ irreducible factors
of degree $2$ and then $\mu(P(X))=(q^2-1)^{h/2}$.
If $h$ is odd, then $P(X)$ has $(h-3)/2$ irreducible factors
of degree $2$ and one irreducible factor of degree 3.
In that case we have $\mu(P(X))=(q^2-1)^{(h-3)/2}(q^3-1)$.
\end{proof}

\begin{theorem}\label{theo2}
We have the following inequalities
\begin{enumerate}
\renewcommand{\theenumi}{\alph{enumi}}
\item
If $q=n-1$ is a prime power, then
$$c(n,h)\leq \frac{\mu(q,q,h+1)}{q-1}.$$
\item If $n-1$ is not a prime power, then
$$
c(n,h)\leq \mu(q,n,h)
$$
where $q$ is the smallest prime power $\geq n$.
\end{enumerate}
\end{theorem}
\begin{proof}
Part (a)  follows from Theorem~\ref{theo1}(b)
and part (b) follows from Theorem~\ref{theo1}(a).
\end{proof}
A lower bound for $A(n,2h+2,w)$ follows from Corollary~\ref{cor1},
Theorem~\ref{theo2} and Lemma~\ref{lem51}. This bound
is almost always superior to the Gilbert-Varshamov type bound for constant
weight codes:
\begin{equation}\label{eq5}
A(n,2h+2,w)\geq \frac{{n\choose w}}{\sum_{i=0}^h {w\choose i}{n-w\choose i}}.
\end{equation}
If $C_w\subseteq \F_2^n$ is a code of constant weight $w$
and minimum distance $\geq 2h+2$ for all $w$, and $u$ is an integer, then
$$
\bigcup_{w\equiv u(\bmod 2h+2)} C_w
$$
is also code with minimum distance $2h+2$. This shows that
\begin{equation}
A(n,2h+2)\geq \sum_{w\equiv u(\bmod 2h+2)}A(n,2h+2,w).
\end{equation}
Heuristically, the best choice for $u$ is $\lfloor n/2\rfloor$
because
$$
\sum_{w\equiv u (m\bmod 2h+2)}{n\choose w}
$$
is maximal for $u=\lfloor n/2\rfloor$.

Using our lower bounds for $A(n,2h+2,w)$ (or the Gilbert-Varshamov
lower bound in those instances where it is better), one obtains
a lower bound for $A(n,2h+2)$. For codes with odd minimum distance
we note the well-known fact that
$$
A(n,2h+1)=A(n+1,2h+2).
$$
\section{Asymptotics}\label{sec4}
\noindent The sphere packing bound tells us $A(n,2h+1)\leq B(n,h)$
where
$$
B(n,h)=\frac{2^n}{{n\choose 0}+{n\choose 1}+\cdots {n\choose h}}=
\frac{2^n h!}{n^h}(1+o(1))
$$
(as $n\to \infty$).
We study the limit densities for $h$-error correcting codes
$$
\rho_{{\rm inf}}(h)=\liminf_{n\to \infty}\frac{A(n,2h+1)}{B(n,h)}
\mbox{ and }
\rho_{{\rm sup}}(h)=\limsup_{n\to \infty}\frac{A(n,2h+1)}{B(n,h)}.
$$
The best known estimates for $\rho_{{\rm inf}}(h)$
and $\rho_{{\rm sup}}(h)$ up to now came from BCH codes:
A BCH code of length $n=2^m-1$ of designed distance $2h+1$
has minimum distance $\geq 2h+1$ and dimension $\geq n-hm$,
so 
$$A(2^m-1,2h+1)\geq 2^{n-hm}=\frac{2^n}{(n+1)^h}
$$
and it follows that 
$$\rho_{{\rm sup}}(h)\geq\frac{1}{h!}.
$$
If $2^{m-1}\leq n \leq 2^{m}-1$ then shortening
the BCH code of length $2^m-1$ above gives a code
of length $n$, minimum distance $\geq 2h+1$ of dimension
$\geq n-hm$.
From this follows that
$$
A(n,2h+1)\geq 2^{n-hm}\geq\frac{2^n}{(2n)^h}
$$
and it follows that
$$
\rho_{{\rm inf}}(h)\geq \frac{1}{h!2^h}.
$$

The bound for $\rho_{{\rm inf}}(h)$ can be improved as follows.
From  Corollary~\ref{cor1},
Theorem~\ref{theo2} and Lemma~\ref{lem51} follows that
$$
A(n,2h+1)=A(n+1,2h+2)\geq \sum_{w\equiv u(\bmod\ 2h+2)}
A(n+1,2h+2,w)\geq 
$$
$$
\geq
\frac{1}{c(n+1,2h+2)}\sum_{w\equiv u(\bmod 2h+2)}
{n+1\choose w}=\frac{1}{c(n+1,2h+2)}\frac{2^{n+1}}{2h+2}(1+o(1))=
$$
$$=
\frac{2^{n}}{(h+1)n^h}(1+o(1))
$$
(where $u=\lfloor (n+1)/2\rfloor$). 
The following theorem follows.
\begin{theorem}
We have the following lower bound for $\rho_{\rm inf}(h)$:
$$
\rho_{{\rm inf}}(h)\geq \frac{1}{(h+1)!}.
$$
\end{theorem}

\section{Table of Improved bounds for $A(n,d)$}\label{sec6}
The following table gives all the improved bounds
for $A(n,d)$ where $1\leq n\leq 512$ and $3\leq d\leq 29$.
The ``new'' column and ``old'' column give the
new and old lower bounds for $\log_2 A(n,d)$ respectively.
The ``ratio'' column gives the ratio of the new and old
lower bounds for $A(n,d)$.

\vfill\pagebreak
{\tiny
\begin{tabular}{||c|c|c|c|c||c|c|c|c|c||c|c|c|c|c||}
\hline\hline
 n & d & new & old & ratio &
n & d & new & old & ratio & n & d & new & old & ratio \\
\hline\hline
279 & 5 & 261.1513 & 261.0000 & 1.1106  & 
280 & 5 & 262.1462 & 262.0000 & 1.1067  & 
281 & 5 & 263.1410 & 263.0000 & 1.1027 \\ 
282 & 5 & 264.1257 & 264.0000 & 1.0911  & 
283 & 5 & 265.1206 & 265.0000 & 1.0872  & 
284 & 5 & 266.0752 & 266.0000 & 1.0535 \\ 
285 & 5 & 267.0752 & 267.0000 & 1.0535  & 
286 & 5 & 268.0752 & 268.0000 & 1.0535  & 
287 & 5 & 269.0702 & 269.0000 & 1.0499 \\ 
288 & 5 & 270.0652 & 270.0000 & 1.0462  & 
289 & 5 & 271.0602 & 271.0000 & 1.0426  & 
290 & 5 & 272.0354 & 272.0000 & 1.0248 \\ 
291 & 5 & 273.0305 & 273.0000 & 1.0213  & 
292 & 5 & 274.0255 & 274.0000 & 1.0179  & 
293 & 5 & 275.0206 & 275.0000 & 1.0144 \\ 
 \hline 312 & 7 & 285.1299 & 285.0000 & 1.0943  & 
313 & 7 & 286.1254 & 286.0000 & 1.0908  & 
314 & 7 & 287.0841 & 287.0000 & 1.0600 \\ 
315 & 7 & 288.0796 & 288.0000 & 1.0567  & 
316 & 7 & 289.0750 & 289.0000 & 1.0533  & 
317 & 7 & 290.0705 & 290.0000 & 1.0500 \\ 
168 & 9 & 136.0752 & 136.0000 & 1.0535  & 
169 & 9 & 137.0666 & 137.0000 & 1.0473  & 
276 & 9 & 241.2231 & 241.0000 & 1.1673 \\ 
277 & 9 & 242.2179 & 242.0000 & 1.1630  & 
278 & 9 & 243.1507 & 243.0000 & 1.1101  & 
279 & 9 & 244.1455 & 244.0000 & 1.1061 \\ 
280 & 9 & 245.1404 & 245.0000 & 1.1022  & 
281 & 9 & 246.1352 & 246.0000 & 1.0983  & 
282 & 9 & 247.0995 & 247.0000 & 1.0714 \\ 
 \hline 283 & 9 & 248.0944 & 248.0000 & 1.0676  & 
316 & 9 & 280.4447 & 280.0000 & 1.3611  & 
317 & 9 & 281.4402 & 281.0000 & 1.3568 \\ 
318 & 9 & 282.2128 & 282.0000 & 1.1589  & 
319 & 9 & 283.2128 & 283.0000 & 1.1589  & 
320 & 9 & 284.2128 & 284.0000 & 1.1589 \\ 
321 & 9 & 285.2128 & 285.0000 & 1.1589  & 
322 & 9 & 286.2128 & 286.0000 & 1.1589  & 
323 & 9 & 287.2128 & 287.0000 & 1.1589 \\ 
324 & 9 & 288.2128 & 288.0000 & 1.1589  & 
325 & 9 & 289.2128 & 289.0000 & 1.1589  & 
326 & 9 & 290.2128 & 290.0000 & 1.1589 \\ 
327 & 9 & 291.2084 & 291.0000 & 1.1554  & 
328 & 9 & 292.2041 & 292.0000 & 1.1519  & 
329 & 9 & 293.1997 & 293.0000 & 1.1485 \\ 
 \hline 330 & 9 & 294.1953 & 294.0000 & 1.1450  & 
331 & 9 & 295.1910 & 295.0000 & 1.1415  & 
332 & 9 & 296.1088 & 296.0000 & 1.0783 \\ 
333 & 9 & 297.1045 & 297.0000 & 1.0751  & 
334 & 9 & 298.1002 & 298.0000 & 1.0720  & 
335 & 9 & 299.0960 & 299.0000 & 1.0688 \\ 
336 & 9 & 300.0917 & 300.0000 & 1.0656  & 
337 & 9 & 301.0874 & 301.0000 & 1.0624  & 
338 & 9 & 302.0067 & 302.0000 & 1.0046 \\ 
339 & 9 & 303.0024 & 303.0000 & 1.0017  & 
150 & 11 & 111.2378 & 111.1699 & 1.0482  & 
151 & 11 & 112.2284 & 112.0000 & 1.1715 \\ 
281 & 11 & 237.7380 & 237.2384 & 1.4138  & 
282 & 11 & 238.6919 & 238.0000 & 1.6154  & 
283 & 11 & 239.6868 & 239.0000 & 1.6098 \\ 
 \hline 284 & 11 & 240.5606 & 240.0000 & 1.4748  & 
285 & 11 & 241.5556 & 241.0000 & 1.4698  & 
286 & 11 & 242.5506 & 242.0000 & 1.4647 \\ 
287 & 11 & 243.5456 & 243.0000 & 1.4596  & 
288 & 11 & 244.5406 & 244.0000 & 1.4545  & 
289 & 11 & 245.5356 & 245.0000 & 1.4495 \\ 
290 & 11 & 246.4513 & 246.0000 & 1.3672  & 
291 & 11 & 247.4463 & 247.0000 & 1.3626  & 
292 & 11 & 248.4414 & 248.0000 & 1.3579 \\ 
293 & 11 & 249.4365 & 249.0000 & 1.3533  & 
294 & 11 & 250.1282 & 250.0000 & 1.0929  & 
295 & 11 & 251.1282 & 251.0000 & 1.0929 \\ 
296 & 11 & 252.1282 & 252.0000 & 1.0929  & 
297 & 11 & 253.1282 & 253.0000 & 1.0929  & 
298 & 11 & 254.1282 & 254.0000 & 1.0929 \\ 
 \hline 299 & 11 & 255.1282 & 255.0000 & 1.0929  & 
300 & 11 & 256.1282 & 256.0000 & 1.0929  & 
301 & 11 & 257.1282 & 257.0000 & 1.0929 \\ 
302 & 11 & 258.1235 & 258.0000 & 1.0894  & 
303 & 11 & 259.1188 & 259.0000 & 1.0858  & 
304 & 11 & 260.1141 & 260.0000 & 1.0823 \\ 
305 & 11 & 261.1094 & 261.0000 & 1.0788  & 
306 & 11 & 262.1046 & 262.0000 & 1.0752  & 
307 & 11 & 263.1000 & 263.0000 & 1.0718 \\ 
308 & 11 & 264.0206 & 264.0000 & 1.0143  & 
309 & 11 & 265.0159 & 265.0000 & 1.0111  & 
310 & 11 & 266.0113 & 266.0000 & 1.0078 \\ 
311 & 11 & 267.0066 & 267.0000 & 1.0046  & 
312 & 11 & 267.9650 & 267.0000 & 1.9521  & 
313 & 11 & 268.9604 & 268.0000 & 1.9459 \\ 
 \hline 314 & 11 & 269.8825 & 269.0000 & 1.8436  & 
315 & 11 & 270.8780 & 270.0000 & 1.8378  & 
316 & 11 & 271.8734 & 271.0000 & 1.8320 \\ 
317 & 11 & 272.8689 & 272.0000 & 1.8262  & 
318 & 11 & 273.5835 & 273.0000 & 1.4984  & 
319 & 11 & 274.5835 & 274.0000 & 1.4984 \\ 
320 & 11 & 275.5835 & 275.0000 & 1.4984  & 
321 & 11 & 276.5835 & 276.0000 & 1.4984  & 
322 & 11 & 277.5835 & 277.0000 & 1.4984 \\ 
323 & 11 & 278.5835 & 278.0000 & 1.4984  & 
324 & 11 & 279.5835 & 279.0000 & 1.4984  & 
325 & 11 & 280.5835 & 280.0000 & 1.4984 \\ 
326 & 11 & 281.5791 & 281.0000 & 1.4939  & 
327 & 11 & 282.5747 & 282.0000 & 1.4894  & 
328 & 11 & 283.5704 & 283.0000 & 1.4849 \\ 
 \hline 329 & 11 & 284.5660 & 284.0000 & 1.4804  & 
330 & 11 & 285.5616 & 285.0000 & 1.4759  & 
331 & 11 & 286.5573 & 286.0000 & 1.4715 \\ 
332 & 11 & 287.4492 & 287.0000 & 1.3653  & 
333 & 11 & 288.4449 & 288.0000 & 1.3612  & 
334 & 11 & 289.4406 & 289.0000 & 1.3572 \\ 
335 & 11 & 290.4364 & 290.0000 & 1.3532  & 
336 & 11 & 291.4320 & 291.0000 & 1.3491  & 
337 & 11 & 292.4278 & 292.0000 & 1.3452 \\ 
338 & 11 & 293.3216 & 293.0000 & 1.2497  & 
339 & 11 & 294.3174 & 294.0000 & 1.2461  & 
340 & 11 & 295.3132 & 295.0000 & 1.2424 \\ 
341 & 11 & 296.3090 & 296.0000 & 1.2388  & 
342 & 11 & 297.3047 & 297.0000 & 1.2352  & 
343 & 11 & 298.3006 & 298.0000 & 1.2316 \\ 
 \hline 344 & 11 & 299.2294 & 299.0000 & 1.1724  & 
345 & 11 & 300.2253 & 300.0000 & 1.1690  & 
346 & 11 & 301.2211 & 301.0000 & 1.1656 \\ 
347 & 11 & 302.2170 & 302.0000 & 1.1623  & 
348 & 11 & 303.1796 & 303.0000 & 1.1326  & 
349 & 11 & 304.1755 & 304.0000 & 1.1294 \\ 
350 & 11 & 305.1056 & 305.0000 & 1.0760  & 
351 & 11 & 306.1015 & 306.0000 & 1.0729  & 
352 & 11 & 307.0974 & 307.0000 & 1.0699 \\ 
353 & 11 & 308.0934 & 308.0000 & 1.0669  & 
157 & 13 & 110.4671 & 110.1699 & 1.2288  & 
158 & 13 & 111.1846 & 111.0000 & 1.1365 \\ 
159 & 13 & 112.1758 & 112.0000 & 1.1296  & 
160 & 13 & 113.1646 & 113.0000 & 1.1208  & 
161 & 13 & 114.1557 & 114.0000 & 1.1140 \\ 
 \hline 162 & 13 & 115.1446 & 115.0000 & 1.1054  & 
163 & 13 & 116.1357 & 116.0000 & 1.0986  & 
281 & 13 & 229.3832 & 229.2384 & 1.1056 \\ 
282 & 13 & 230.3269 & 230.0000 & 1.2543  & 
283 & 13 & 231.3218 & 231.0000 & 1.2498  & 
284 & 13 & 232.1651 & 232.0000 & 1.1213 \\ 
285 & 13 & 233.1601 & 233.0000 & 1.1174  & 
286 & 13 & 234.1550 & 234.0000 & 1.1135  & 
287 & 13 & 235.1500 & 235.0000 & 1.1096 \\ 
288 & 13 & 236.1450 & 236.0000 & 1.1057  & 
289 & 13 & 237.1400 & 236.2384 & 1.8681  & 
290 & 13 & 238.0357 & 237.0000 & 2.0502 \\ 
291 & 13 & 239.0308 & 238.0000 & 2.0432  & 
292 & 13 & 240.0258 & 239.0000 & 2.0361  & 
293 & 13 & 241.0209 & 240.0000 & 2.0291 \\ 
 \hline 294 & 13 & 241.6499 & 241.0000 & 1.5690  & 
295 & 13 & 242.6499 & 242.0000 & 1.5690  & 
296 & 13 & 243.6498 & 243.0000 & 1.5690 \\ 
297 & 13 & 244.6498 & 244.0000 & 1.5690  & 
298 & 13 & 245.6497 & 245.0000 & 1.5689  & 
299 & 13 & 246.6497 & 246.0000 & 1.5689 \\ 
300 & 13 & 247.6497 & 247.0000 & 1.5688  & 
301 & 13 & 248.6450 & 248.0000 & 1.5637  & 
302 & 13 & 249.6402 & 249.0000 & 1.5585 \\ 
303 & 13 & 250.6355 & 250.0000 & 1.5535  & 
304 & 13 & 251.6307 & 251.0000 & 1.5484  & 
305 & 13 & 252.6260 & 252.0000 & 1.5433 \\ 
306 & 13 & 253.6213 & 253.0000 & 1.5382  & 
307 & 13 & 254.6166 & 254.0000 & 1.5332  & 
308 & 13 & 255.5184 & 255.0000 & 1.4324 \\ 
 \hline 309 & 13 & 256.5138 & 256.0000 & 1.4278  & 
310 & 13 & 257.5091 & 257.0000 & 1.4232  & 
311 & 13 & 258.5044 & 258.0000 & 1.4186 \\ 
312 & 13 & 259.4536 & 259.0000 & 1.3694  & 
313 & 13 & 260.4489 & 260.0000 & 1.3650  & 
314 & 13 & 261.3527 & 261.0000 & 1.2769 \\ 
315 & 13 & 262.3481 & 262.0000 & 1.2729  & 
316 & 13 & 263.3436 & 263.0000 & 1.2689  & 
317 & 13 & 264.3390 & 264.0000 & 1.2649 \\ 
328 & 13 & 274.9779 & 274.0000 & 1.9697  & 
329 & 13 & 275.9736 & 275.0000 & 1.9637  & 
330 & 13 & 276.9692 & 276.0000 & 1.9577 \\ 
331 & 13 & 277.9648 & 277.0000 & 1.9518  & 
332 & 13 & 278.8308 & 278.0000 & 1.7786  & 
333 & 13 & 279.8265 & 279.0000 & 1.7734 \\ 
 \hline \hline 
\end{tabular}

\begin{tabular}{||c|c|c|c|c||c|c|c|c|c||c|c|c|c|c||}
\hline\hline
 n & d & new & old & ratio &
n & d & new & old & ratio & n & d & new & old & ratio \\
\hline\hline
334 & 13 & 280.8222 & 280.0000 & 1.7681  & 
335 & 13 & 281.8179 & 281.0000 & 1.7628  & 
336 & 13 & 282.8136 & 282.0000 & 1.7576 \\ 
337 & 13 & 283.8093 & 283.0000 & 1.7524  & 
338 & 13 & 284.6776 & 284.0000 & 1.5995  & 
339 & 13 & 285.6734 & 285.0000 & 1.5948 \\ 
340 & 13 & 286.6692 & 286.0000 & 1.5902  & 
341 & 13 & 287.6650 & 287.0000 & 1.5855  & 
342 & 13 & 288.6608 & 288.0000 & 1.5809 \\ 
343 & 13 & 289.6565 & 289.0000 & 1.5763  & 
344 & 13 & 290.5687 & 290.0000 & 1.4832  & 
345 & 13 & 291.5645 & 291.0000 & 1.4789 \\ 
346 & 13 & 292.5603 & 292.0000 & 1.4746  & 
347 & 13 & 293.5562 & 293.0000 & 1.4704  & 
348 & 13 & 294.5106 & 294.0000 & 1.4246 \\ 
 \hline 349 & 13 & 295.5064 & 295.0000 & 1.4205  & 
350 & 13 & 296.4201 & 296.0000 & 1.3380  & 
351 & 13 & 297.4160 & 297.0000 & 1.3342 \\ 
352 & 13 & 298.4119 & 298.0000 & 1.3304  & 
353 & 13 & 299.4078 & 299.0000 & 1.3267  & 
354 & 13 & 300.2821 & 300.0000 & 1.2159 \\ 
355 & 13 & 301.2780 & 301.0000 & 1.2125  & 
356 & 13 & 302.2740 & 302.0000 & 1.2092  & 
357 & 13 & 303.2700 & 303.0000 & 1.2058 \\ 
358 & 13 & 304.2659 & 304.0000 & 1.2024  & 
359 & 13 & 305.2619 & 305.0000 & 1.1991  & 
360 & 13 & 306.2178 & 306.0000 & 1.1630 \\ 
361 & 13 & 307.2138 & 307.0000 & 1.1598  & 
362 & 13 & 308.0909 & 308.0000 & 1.0650  & 
363 & 13 & 309.0869 & 309.0000 & 1.0621 \\ 
 \hline 364 & 13 & 310.0830 & 310.0000 & 1.0592  & 
365 & 13 & 311.0791 & 311.0000 & 1.0563  & 
366 & 13 & 312.0751 & 312.0000 & 1.0534 \\ 
367 & 13 & 313.0712 & 313.0000 & 1.0506  & 
164 & 15 & 109.4423 & 109.1699 & 1.2078  & 
165 & 15 & 110.4337 & 110.0000 & 1.3507 \\ 
166 & 15 & 111.4209 & 111.0000 & 1.3388  & 
167 & 15 & 112.4124 & 112.0000 & 1.3309  & 
168 & 15 & 113.2968 & 113.0000 & 1.2284 \\ 
169 & 15 & 114.2883 & 114.0000 & 1.2212  & 
170 & 15 & 115.0735 & 115.0000 & 1.0523  & 
171 & 15 & 116.0652 & 116.0000 & 1.0462 \\ 
172 & 15 & 117.0531 & 117.0000 & 1.0375  & 
173 & 15 & 118.0449 & 118.0000 & 1.0316  & 
289 & 15 & 228.7809 & 228.2384 & 1.4565 \\ 
 \hline 290 & 15 & 229.6566 & 229.0000 & 1.5763  & 
291 & 15 & 230.6516 & 230.0000 & 1.5710  & 
292 & 15 & 231.6463 & 231.0000 & 1.5652 \\ 
293 & 15 & 232.6414 & 232.0000 & 1.5599  & 
294 & 15 & 233.2075 & 233.0000 & 1.1547  & 
295 & 15 & 234.2075 & 234.0000 & 1.1547 \\ 
296 & 15 & 235.2072 & 235.0000 & 1.1544  & 
297 & 15 & 236.2072 & 235.2384 & 1.9571  & 
298 & 15 & 237.2068 & 236.0000 & 2.3083 \\ 
299 & 15 & 238.2068 & 237.0000 & 2.3083  & 
300 & 15 & 239.2018 & 238.0000 & 2.3002  & 
301 & 15 & 240.1971 & 239.0000 & 2.2928 \\ 
302 & 15 & 241.1921 & 240.0000 & 2.2848  & 
303 & 15 & 242.1874 & 241.0000 & 2.2774  & 
304 & 15 & 243.1824 & 242.0000 & 2.2695 \\ 
 \hline 305 & 15 & 244.1777 & 243.0000 & 2.2621  & 
306 & 15 & 245.1727 & 244.0000 & 2.2543  & 
307 & 15 & 246.1680 & 245.0000 & 2.2470 \\ 
308 & 15 & 247.0510 & 246.0000 & 2.0719  & 
309 & 15 & 248.0463 & 247.0000 & 2.0653  & 
310 & 15 & 249.0414 & 248.0000 & 2.0582 \\ 
311 & 15 & 250.0368 & 249.0000 & 2.0516  & 
312 & 15 & 250.9764 & 250.0000 & 1.9676  & 
313 & 15 & 251.9718 & 251.0000 & 1.9613 \\ 
314 & 15 & 252.8570 & 252.0000 & 1.8113  & 
315 & 15 & 253.8525 & 253.0000 & 1.8056  & 
316 & 15 & 254.8477 & 254.0000 & 1.7996 \\ 
317 & 15 & 255.8431 & 255.0000 & 1.7940  & 
318 & 15 & 256.4415 & 256.0000 & 1.3581  & 
319 & 15 & 257.4415 & 257.0000 & 1.3581 \\ 
 \hline 320 & 15 & 258.4413 & 258.0000 & 1.3578  & 
321 & 15 & 259.4413 & 259.0000 & 1.3578  & 
322 & 15 & 260.4411 & 260.0000 & 1.3576 \\ 
323 & 15 & 261.4411 & 261.0000 & 1.3576  & 
324 & 15 & 262.4365 & 262.0000 & 1.3534  & 
325 & 15 & 263.4322 & 263.0000 & 1.3493 \\ 
326 & 15 & 264.4276 & 264.0000 & 1.3450  & 
327 & 15 & 265.4233 & 265.0000 & 1.3410  & 
328 & 15 & 266.4187 & 265.0000 & 2.6735 \\ 
329 & 15 & 267.4144 & 266.0000 & 2.6654  & 
330 & 15 & 268.4098 & 267.0000 & 2.6570  & 
331 & 15 & 269.4055 & 268.0000 & 2.6490 \\ 
332 & 15 & 270.2453 & 269.0000 & 2.3708  & 
333 & 15 & 271.2411 & 270.0000 & 2.3637  & 
334 & 15 & 272.2366 & 271.0000 & 2.3565 \\ 
 \hline 335 & 15 & 273.2323 & 272.0000 & 2.3495  & 
336 & 15 & 274.2279 & 273.0000 & 2.3422  & 
337 & 15 & 275.2236 & 274.0000 & 2.3353 \\ 
338 & 15 & 276.0663 & 275.0000 & 2.0941  & 
339 & 15 & 277.0621 & 276.0000 & 2.0880  & 
340 & 15 & 278.0578 & 277.0000 & 2.0817 \\ 
341 & 15 & 279.0536 & 278.0000 & 2.0757  & 
342 & 15 & 280.0492 & 279.0000 & 2.0694  & 
343 & 15 & 281.0450 & 280.0000 & 2.0634 \\ 
344 & 15 & 281.9403 & 281.0000 & 1.9189  & 
345 & 15 & 282.9362 & 282.0000 & 1.9134  & 
346 & 15 & 283.9318 & 283.0000 & 1.9077 \\ 
347 & 15 & 284.9277 & 284.0000 & 1.9023  & 
348 & 15 & 285.8737 & 285.0000 & 1.8323  & 
349 & 15 & 286.8696 & 286.0000 & 1.8271 \\ 
 \hline 350 & 15 & 287.7667 & 287.0000 & 1.7013  & 
351 & 15 & 288.7626 & 288.0000 & 1.6965  & 
352 & 15 & 289.7584 & 289.0000 & 1.6916 \\ 
353 & 15 & 290.7543 & 290.0000 & 1.6868  & 
354 & 15 & 291.6041 & 291.0000 & 1.5200  & 
355 & 15 & 292.6001 & 292.0000 & 1.5158 \\ 
356 & 15 & 293.5960 & 293.0000 & 1.5115  & 
357 & 15 & 294.5920 & 294.0000 & 1.5073  & 
358 & 15 & 295.5878 & 295.0000 & 1.5030 \\ 
359 & 15 & 296.5838 & 296.0000 & 1.4988  & 
360 & 15 & 297.5316 & 297.0000 & 1.4455  & 
361 & 15 & 298.5276 & 298.0000 & 1.4416 \\ 
362 & 15 & 299.3808 & 299.0000 & 1.3020  & 
363 & 15 & 300.3768 & 300.0000 & 1.2985  & 
364 & 15 & 301.3728 & 301.0000 & 1.2949 \\ 
 \hline 365 & 15 & 302.3689 & 302.0000 & 1.2914  & 
366 & 15 & 303.3649 & 303.0000 & 1.2878  & 
367 & 15 & 304.3609 & 304.0000 & 1.2843 \\ 
368 & 15 & 305.2165 & 305.0000 & 1.1619  & 
369 & 15 & 306.2126 & 306.0000 & 1.1588  & 
370 & 15 & 307.2087 & 307.0000 & 1.1556 \\ 
371 & 15 & 308.2048 & 308.0000 & 1.1525  & 
372 & 15 & 309.2008 & 309.0000 & 1.1494  & 
373 & 15 & 310.1970 & 310.0000 & 1.1463 \\ 
374 & 15 & 311.0548 & 311.0000 & 1.0387  & 
375 & 15 & 312.0510 & 312.0000 & 1.0360  & 
376 & 15 & 313.0471 & 313.0000 & 1.0332 \\ 
377 & 15 & 314.0433 & 314.0000 & 1.0305  & 
378 & 15 & 315.0395 & 315.0000 & 1.0277  & 
379 & 15 & 316.0357 & 316.0000 & 1.0250 \\ 
 \hline 171 & 17 & 108.5553 & 108.1699 & 1.3062  & 
172 & 17 & 109.5415 & 109.0000 & 1.4555  & 
173 & 17 & 110.5331 & 110.0000 & 1.4471 \\ 
174 & 17 & 111.1749 & 111.0000 & 1.1289  & 
175 & 17 & 112.1668 & 112.0000 & 1.1226  & 
176 & 17 & 113.1536 & 113.0000 & 1.1123 \\ 
177 & 17 & 114.1455 & 114.0000 & 1.1061  & 
178 & 17 & 115.1324 & 115.0000 & 1.0961  & 
179 & 17 & 116.1243 & 116.0000 & 1.0900 \\ 
192 & 17 & 128.2315 & 128.0000 & 1.1741  & 
193 & 17 & 129.2240 & 129.0000 & 1.1680  & 
194 & 17 & 130.0053 & 130.0000 & 1.0037 \\ 
289 & 17 & 220.4594 & 220.2384 & 1.1655  & 
290 & 17 & 221.3146 & 221.0000 & 1.2437  & 
291 & 17 & 222.3097 & 222.0000 & 1.2394 \\ 
 \hline 292 & 17 & 223.3038 & 223.0000 & 1.2344  & 
293 & 17 & 224.2988 & 224.0000 & 1.2301  & 
297 & 17 & 227.8008 & 227.2384 & 1.4767 \\ 
298 & 17 & 228.7999 & 228.0000 & 1.7410  & 
299 & 17 & 229.7952 & 229.0000 & 1.7353  & 
300 & 17 & 230.7896 & 230.0000 & 1.7286 \\ 
301 & 17 & 231.7849 & 231.0000 & 1.7230  & 
302 & 17 & 232.7794 & 232.0000 & 1.7164  & 
303 & 17 & 233.7747 & 233.0000 & 1.7108 \\ 
304 & 17 & 234.7692 & 234.0000 & 1.7043  & 
305 & 17 & 235.7645 & 234.2384 & 2.8800  & 
306 & 17 & 236.7590 & 235.0000 & 3.3846 \\ 
307 & 17 & 237.7543 & 236.0000 & 3.3736  & 
308 & 17 & 238.6181 & 237.0000 & 3.0697  & 
309 & 17 & 239.6135 & 238.0000 & 3.0598 \\ 
 \hline 310 & 17 & 240.6081 & 239.0000 & 3.0485  & 
311 & 17 & 241.6034 & 240.0000 & 3.0386  & 
312 & 17 & 242.5334 & 241.0000 & 2.8946 \\ 
313 & 17 & 243.5287 & 242.0000 & 2.8853  & 
314 & 17 & 244.3952 & 243.0000 & 2.6303  & 
315 & 17 & 245.3907 & 244.0000 & 2.6220 \\ 
316 & 17 & 246.3854 & 245.0000 & 2.6125  & 
317 & 17 & 247.3809 & 246.0000 & 2.6042  & 
318 & 17 & 247.9209 & 247.0000 & 1.8932 \\ 
319 & 17 & 248.9209 & 248.0000 & 1.8932  & 
320 & 17 & 249.9202 & 249.0000 & 1.8924  & 
321 & 17 & 250.9202 & 250.0000 & 1.8924 \\ 
322 & 17 & 251.9196 & 251.0000 & 1.8916  & 
323 & 17 & 252.9152 & 252.0000 & 1.8859  & 
324 & 17 & 253.9103 & 253.0000 & 1.8794 \\ 
 \hline \hline 
\end{tabular}

\begin{tabular}{||c|c|c|c|c||c|c|c|c|c||c|c|c|c|c||}
\hline\hline
 n & d & new & old & ratio &
n & d & new & old & ratio & n & d & new & old & ratio \\
\hline\hline
325 & 17 & 254.9059 & 254.0000 & 1.8737  & 
326 & 17 & 255.9010 & 255.0000 & 1.8673  & 
327 & 17 & 256.8966 & 256.0000 & 1.8617 \\ 
328 & 17 & 257.8917 & 257.0000 & 1.8553  & 
329 & 17 & 258.8873 & 258.0000 & 1.8497  & 
330 & 17 & 259.8824 & 259.0000 & 1.8435 \\ 
331 & 17 & 260.8780 & 260.0000 & 1.8379  & 
332 & 17 & 261.6917 & 260.0000 & 3.2303  & 
333 & 17 & 262.6874 & 261.0000 & 3.2207 \\ 
334 & 17 & 263.6826 & 262.0000 & 3.2100  & 
335 & 17 & 264.6783 & 263.0000 & 3.2005  & 
336 & 17 & 265.6735 & 264.0000 & 3.1899 \\ 
337 & 17 & 266.6692 & 265.0000 & 3.1804  & 
338 & 17 & 267.4862 & 266.0000 & 2.8015  & 
339 & 17 & 268.4820 & 267.0000 & 2.7933 \\ 
 \hline 340 & 17 & 269.4773 & 268.0000 & 2.7843  & 
341 & 17 & 270.4731 & 269.0000 & 2.7761  & 
342 & 17 & 271.4684 & 270.0000 & 2.7672 \\ 
343 & 17 & 272.4642 & 271.0000 & 2.7591  & 
344 & 17 & 273.3425 & 272.0000 & 2.5359  & 
345 & 17 & 274.3383 & 273.0000 & 2.5286 \\ 
346 & 17 & 275.3337 & 274.0000 & 2.5205  & 
347 & 17 & 276.3296 & 275.0000 & 2.5133  & 
348 & 17 & 277.2670 & 276.0000 & 2.4066 \\ 
349 & 17 & 278.2628 & 277.0000 & 2.3997  & 
350 & 17 & 279.1432 & 278.0000 & 2.2087  & 
351 & 17 & 280.1391 & 279.0000 & 2.2025 \\ 
352 & 17 & 281.1347 & 280.0000 & 2.1957  & 
353 & 17 & 282.1305 & 281.0000 & 2.1894  & 
354 & 17 & 282.9558 & 282.0000 & 1.9397 \\ 
 \hline 355 & 17 & 283.9518 & 283.0000 & 1.9343  & 
356 & 17 & 284.9474 & 284.0000 & 1.9284  & 
357 & 17 & 285.9434 & 285.0000 & 1.9230 \\ 
358 & 17 & 286.9390 & 286.0000 & 1.9172  & 
359 & 17 & 287.9350 & 287.0000 & 1.9119  & 
360 & 17 & 288.8745 & 288.0000 & 1.8334 \\ 
361 & 17 & 289.8705 & 289.0000 & 1.8283  & 
362 & 17 & 290.6997 & 290.0000 & 1.6241  & 
363 & 17 & 291.6957 & 291.0000 & 1.6197 \\ 
364 & 17 & 292.6915 & 292.0000 & 1.6149  & 
365 & 17 & 293.6876 & 293.0000 & 1.6105  & 
366 & 17 & 294.6833 & 294.0000 & 1.6058 \\ 
367 & 17 & 295.6794 & 295.0000 & 1.6014  & 
368 & 17 & 296.5113 & 296.0000 & 1.4253  & 
369 & 17 & 297.5074 & 297.0000 & 1.4215 \\ 
 \hline 370 & 17 & 298.5033 & 298.0000 & 1.4174  & 
371 & 17 & 299.4994 & 299.0000 & 1.4136  & 
372 & 17 & 300.4952 & 300.0000 & 1.4096 \\ 
373 & 17 & 301.4914 & 301.0000 & 1.4058  & 
374 & 17 & 302.3260 & 302.0000 & 1.2535  & 
375 & 17 & 303.3222 & 303.0000 & 1.2502 \\ 
376 & 17 & 304.3181 & 304.0000 & 1.2467  & 
377 & 17 & 305.3143 & 305.0000 & 1.2434  & 
378 & 17 & 306.3102 & 306.0000 & 1.2399 \\ 
379 & 17 & 307.3064 & 307.0000 & 1.2366  & 
380 & 17 & 308.1963 & 308.0000 & 1.1458  & 
381 & 17 & 309.1926 & 309.0000 & 1.1428 \\ 
382 & 17 & 310.1886 & 310.0000 & 1.1396  & 
383 & 17 & 311.1848 & 311.0000 & 1.1366  & 
384 & 17 & 312.0237 & 312.0000 & 1.0166 \\ 
 \hline 385 & 17 & 313.0200 & 313.0000 & 1.0140  & 
386 & 17 & 314.0161 & 314.0000 & 1.0112  & 
387 & 17 & 315.0124 & 315.0000 & 1.0086 \\ 
388 & 17 & 316.0085 & 316.0000 & 1.0059  & 
389 & 17 & 317.0047 & 317.0000 & 1.0033  & 
290 & 19 & 213.0107 & 213.0000 & 1.0074 \\ 
291 & 19 & 214.0057 & 214.0000 & 1.0040  & 
297 & 19 & 219.4318 & 219.2384 & 1.1434  & 
298 & 19 & 220.4254 & 220.0000 & 1.3429 \\ 
299 & 19 & 221.4207 & 221.0000 & 1.3386  & 
300 & 19 & 222.4144 & 222.0000 & 1.3327  & 
301 & 19 & 223.4097 & 223.0000 & 1.3284 \\ 
302 & 19 & 224.4034 & 224.0000 & 1.3226  & 
303 & 19 & 225.3987 & 225.0000 & 1.3183  & 
304 & 19 & 226.3924 & 226.0000 & 1.3126 \\ 
 \hline 305 & 19 & 227.3877 & 226.2384 & 2.2181  & 
306 & 19 & 228.3815 & 227.0000 & 2.6053  & 
307 & 19 & 229.3768 & 228.0000 & 2.5969 \\ 
308 & 19 & 230.2212 & 229.0000 & 2.3314  & 
309 & 19 & 231.2165 & 230.0000 & 2.3239  & 
310 & 19 & 232.2104 & 231.0000 & 2.3141 \\ 
311 & 19 & 233.2058 & 232.0000 & 2.3067  & 
312 & 19 & 234.1258 & 233.0000 & 2.1822  & 
313 & 19 & 235.1212 & 234.0000 & 2.1753 \\ 
314 & 19 & 235.9686 & 235.0000 & 1.9570  & 
315 & 19 & 236.9641 & 236.0000 & 1.9508  & 
316 & 19 & 237.9581 & 237.0000 & 1.9428 \\ 
317 & 19 & 238.9536 & 238.0000 & 1.9367  & 
318 & 19 & 239.4349 & 239.0000 & 1.3518  & 
319 & 19 & 240.4349 & 240.0000 & 1.3518 \\ 
 \hline 320 & 19 & 241.4336 & 240.0000 & 2.7013  & 
321 & 19 & 242.4336 & 241.0000 & 2.7013  & 
322 & 19 & 243.4280 & 242.0000 & 2.6908 \\ 
323 & 19 & 244.4237 & 243.0000 & 2.6827  & 
324 & 19 & 245.4181 & 244.0000 & 2.6723  & 
325 & 19 & 246.4137 & 245.0000 & 2.6642 \\ 
326 & 19 & 247.4081 & 246.0000 & 2.6540  & 
327 & 19 & 248.4038 & 247.0000 & 2.6460  & 
328 & 19 & 249.3983 & 248.0000 & 2.6358 \\ 
329 & 19 & 250.3939 & 249.0000 & 2.6279  & 
330 & 19 & 251.3884 & 250.0000 & 2.6178  & 
331 & 19 & 252.3840 & 251.0000 & 2.6100 \\ 
332 & 19 & 253.1711 & 252.0000 & 2.2519  & 
333 & 19 & 254.1669 & 253.0000 & 2.2452  & 
334 & 19 & 255.1615 & 254.0000 & 2.2369 \\ 
 \hline 335 & 19 & 256.1572 & 255.0000 & 2.2303  & 
336 & 19 & 257.1518 & 256.0000 & 2.2220  & 
337 & 19 & 258.1476 & 257.0000 & 2.2154 \\ 
338 & 19 & 258.9385 & 258.0000 & 1.9165  & 
339 & 19 & 259.9343 & 259.0000 & 1.9110  & 
340 & 19 & 260.9291 & 260.0000 & 1.9040 \\ 
341 & 19 & 261.9249 & 261.0000 & 1.8985  & 
342 & 19 & 262.9197 & 262.0000 & 1.8917  & 
343 & 19 & 263.9155 & 263.0000 & 1.8862 \\ 
344 & 19 & 264.7765 & 264.0000 & 1.7129  & 
345 & 19 & 265.7723 & 265.0000 & 1.7080  & 
346 & 19 & 266.7672 & 266.0000 & 1.7020 \\ 
347 & 19 & 267.7631 & 267.0000 & 1.6971  & 
348 & 19 & 268.6916 & 268.0000 & 1.6151  & 
349 & 19 & 269.6875 & 269.0000 & 1.6105 \\ 
 \hline 350 & 19 & 270.5510 & 270.0000 & 1.4651  & 
351 & 19 & 271.5469 & 271.0000 & 1.4609  & 
352 & 19 & 272.5419 & 272.0000 & 1.4559 \\ 
353 & 19 & 273.5378 & 273.0000 & 1.4518  & 
354 & 19 & 274.3383 & 274.0000 & 1.2642  & 
355 & 19 & 275.3343 & 275.0000 & 1.2607 \\ 
356 & 19 & 276.3294 & 276.0000 & 1.2565  & 
357 & 19 & 277.3254 & 277.0000 & 1.2530  & 
358 & 19 & 278.3205 & 278.0000 & 1.2488 \\ 
359 & 19 & 279.3165 & 279.0000 & 1.2453  & 
360 & 19 & 280.2476 & 280.0000 & 1.1872  & 
361 & 19 & 281.2436 & 281.0000 & 1.1840 \\ 
362 & 19 & 282.0485 & 282.0000 & 1.0342  & 
363 & 19 & 283.0446 & 283.0000 & 1.0314  & 
364 & 19 & 284.0399 & 284.0000 & 1.0280 \\ 
 \hline 365 & 19 & 285.0360 & 285.0000 & 1.0253  & 
366 & 19 & 286.0313 & 286.0000 & 1.0219  & 
367 & 19 & 287.0274 & 287.0000 & 1.0192 \\ 
384 & 19 & 303.2843 & 303.0000 & 1.2178  & 
385 & 19 & 304.2806 & 304.0000 & 1.2147  & 
386 & 19 & 305.2763 & 305.0000 & 1.2111 \\ 
387 & 19 & 306.2726 & 306.0000 & 1.2080  & 
388 & 19 & 307.2683 & 307.0000 & 1.2044  & 
389 & 19 & 308.2646 & 308.0000 & 1.2013 \\ 
390 & 19 & 309.0252 & 309.0000 & 1.0176  & 
391 & 19 & 310.0216 & 310.0000 & 1.0151  & 
392 & 19 & 311.0174 & 311.0000 & 1.0122 \\ 
393 & 19 & 312.0138 & 312.0000 & 1.0096  & 
394 & 19 & 313.0096 & 313.0000 & 1.0067  & 
395 & 19 & 314.0060 & 314.0000 & 1.0042 \\ 
 \hline 396 & 19 & 315.0018 & 315.0000 & 1.0013  & 
298 & 21 & 212.0855 & 212.0000 & 1.0610  & 
299 & 21 & 213.0808 & 213.0000 & 1.0576 \\ 
300 & 21 & 214.0736 & 214.0000 & 1.0523  & 
301 & 21 & 215.0689 & 215.0000 & 1.0489  & 
302 & 21 & 216.0618 & 216.0000 & 1.0437 \\ 
303 & 21 & 217.0571 & 217.0000 & 1.0403  & 
304 & 21 & 218.0500 & 218.0000 & 1.0352  & 
305 & 21 & 219.0453 & 218.2384 & 1.7494 \\ 
306 & 21 & 220.0382 & 219.0000 & 2.0537  & 
307 & 21 & 221.0335 & 220.0000 & 2.0470  & 
308 & 21 & 221.8584 & 221.0000 & 1.8131 \\ 
309 & 21 & 222.8538 & 222.0000 & 1.8072  & 
310 & 21 & 223.8469 & 223.0000 & 1.7986  & 
311 & 21 & 224.8422 & 224.0000 & 1.7928 \\ 
 \hline 312 & 21 & 225.7522 & 225.0000 & 1.6843  & 
313 & 21 & 226.7475 & 225.2384 & 2.8464  & 
314 & 21 & 227.5759 & 226.0000 & 2.9812 \\ 
315 & 21 & 228.5713 & 227.0000 & 2.9718  & 
316 & 21 & 229.5647 & 228.0000 & 2.9581  & 
317 & 21 & 230.5601 & 229.0000 & 2.9487 \\ 
318 & 21 & 230.9827 & 230.0000 & 1.9761  & 
319 & 21 & 231.9827 & 231.0000 & 1.9761  & 
320 & 21 & 232.9806 & 232.0000 & 1.9733 \\ 
321 & 21 & 233.9762 & 232.2384 & 3.3353  & 
322 & 21 & 234.9699 & 233.0000 & 3.9173  & 
323 & 21 & 235.9655 & 234.0000 & 3.9055 \\ 
324 & 21 & 236.9592 & 235.0000 & 3.8884  & 
325 & 21 & 237.9548 & 236.0000 & 3.8766  & 
326 & 21 & 238.9485 & 237.0000 & 3.8597 \\ 
 \hline \hline 
\end{tabular}

\begin{tabular}{||c|c|c|c|c||c|c|c|c|c||c|c|c|c|c||}
\hline\hline
 n & d & new & old & ratio &
n & d & new & old & ratio & n & d & new & old & ratio \\
\hline\hline
327 & 21 & 239.9441 & 238.0000 & 3.8481  & 
328 & 21 & 240.9379 & 239.0000 & 3.8314  & 
329 & 21 & 241.9335 & 240.0000 & 3.8198 \\ 
330 & 21 & 242.9273 & 241.0000 & 3.8034  & 
331 & 21 & 243.9229 & 242.0000 & 3.7919  & 
332 & 21 & 244.6834 & 243.0000 & 3.2118 \\ 
333 & 21 & 245.6791 & 244.0000 & 3.2023  & 
334 & 21 & 246.6730 & 245.0000 & 3.1889  & 
335 & 21 & 247.6688 & 246.0000 & 3.1794 \\ 
336 & 21 & 248.6627 & 247.0000 & 3.1661  & 
337 & 21 & 249.6584 & 248.0000 & 3.1567  & 
338 & 21 & 250.4232 & 249.0000 & 2.6818 \\ 
339 & 21 & 251.4190 & 250.0000 & 2.6740  & 
340 & 21 & 252.4131 & 251.0000 & 2.6631  & 
341 & 21 & 253.4089 & 252.0000 & 2.6553 \\ 
 \hline 342 & 21 & 254.4030 & 253.0000 & 2.6445  & 
343 & 21 & 255.3988 & 254.0000 & 2.6368  & 
344 & 21 & 256.2424 & 255.0000 & 2.3659 \\ 
345 & 21 & 257.2383 & 256.0000 & 2.3591  & 
346 & 21 & 258.2325 & 257.0000 & 2.3497  & 
347 & 21 & 259.2283 & 258.0000 & 2.3429 \\ 
348 & 21 & 260.1480 & 259.0000 & 2.2161  & 
349 & 21 & 261.1439 & 260.0000 & 2.2097  & 
350 & 21 & 261.9902 & 261.0000 & 1.9865 \\ 
351 & 21 & 262.9861 & 262.0000 & 1.9809  & 
352 & 21 & 263.9805 & 263.0000 & 1.9732  & 
353 & 21 & 264.9764 & 264.0000 & 1.9676 \\ 
354 & 21 & 265.7520 & 265.0000 & 1.6841  & 
355 & 21 & 266.7479 & 266.0000 & 1.6794  & 
356 & 21 & 267.7425 & 267.0000 & 1.6730 \\ 
 \hline 357 & 21 & 268.7385 & 268.0000 & 1.6684  & 
358 & 21 & 269.7330 & 269.0000 & 1.6621  & 
359 & 21 & 270.7290 & 270.0000 & 1.6575 \\ 
360 & 21 & 271.6515 & 271.0000 & 1.5708  & 
361 & 21 & 272.6475 & 272.0000 & 1.5664  & 
362 & 21 & 273.4280 & 273.0000 & 1.3454 \\ 
363 & 21 & 274.4241 & 274.0000 & 1.3417  & 
364 & 21 & 275.4188 & 275.0000 & 1.3368  & 
365 & 21 & 276.4149 & 276.0000 & 1.3332 \\ 
366 & 21 & 277.4096 & 277.0000 & 1.3283  & 
367 & 21 & 278.4057 & 278.0000 & 1.3247  & 
368 & 21 & 279.1898 & 279.0000 & 1.1406 \\ 
369 & 21 & 280.1860 & 280.0000 & 1.1376  & 
370 & 21 & 281.1808 & 281.0000 & 1.1335  & 
371 & 21 & 282.1770 & 282.0000 & 1.1305 \\ 
 \hline 372 & 21 & 283.1718 & 283.0000 & 1.1265  & 
373 & 21 & 284.1680 & 284.0000 & 1.1235  & 
384 & 21 & 294.5738 & 294.0000 & 1.4884 \\ 
385 & 21 & 295.5700 & 295.0000 & 1.4846  & 
386 & 21 & 296.5653 & 296.0000 & 1.4796  & 
387 & 21 & 297.5615 & 297.0000 & 1.4758 \\ 
388 & 21 & 298.5568 & 298.0000 & 1.4710  & 
389 & 21 & 299.5530 & 299.0000 & 1.4672  & 
390 & 21 & 300.2838 & 300.0000 & 1.2174 \\ 
391 & 21 & 301.2802 & 301.0000 & 1.2144  & 
392 & 21 & 302.2755 & 302.0000 & 1.2105  & 
393 & 21 & 303.2719 & 303.0000 & 1.2074 \\ 
394 & 21 & 304.2673 & 304.0000 & 1.2035  & 
395 & 21 & 305.2636 & 305.0000 & 1.2005  & 
396 & 21 & 306.2590 & 306.0000 & 1.1967 \\ 
 \hline 397 & 21 & 307.2554 & 307.0000 & 1.1937  & 
398 & 21 & 308.1206 & 308.0000 & 1.0872  & 
399 & 21 & 309.1170 & 309.0000 & 1.0845 \\ 
400 & 21 & 310.1125 & 310.0000 & 1.0811  & 
401 & 21 & 311.1089 & 311.0000 & 1.0784  & 
305 & 23 & 210.7313 & 210.2384 & 1.4072 \\ 
306 & 23 & 211.7235 & 211.0000 & 1.6511  & 
307 & 23 & 212.7188 & 212.0000 & 1.6458  & 
308 & 23 & 213.5242 & 213.0000 & 1.4382 \\ 
309 & 23 & 214.5196 & 214.0000 & 1.4336  & 
310 & 23 & 215.5119 & 215.0000 & 1.4260  & 
311 & 23 & 216.5073 & 216.0000 & 1.4214 \\ 
312 & 23 & 217.4072 & 217.0000 & 1.3262  & 
313 & 23 & 218.4027 & 217.2384 & 2.2412  & 
314 & 23 & 219.2119 & 218.0000 & 2.3164 \\ 
 \hline 315 & 23 & 220.2074 & 219.0000 & 2.3092  & 
316 & 23 & 221.1999 & 220.0000 & 2.2973  & 
317 & 23 & 222.1954 & 221.0000 & 2.2901 \\ 
318 & 23 & 222.5592 & 222.0000 & 1.4735  & 
319 & 23 & 223.5592 & 223.0000 & 1.4735  & 
320 & 23 & 224.5521 & 224.0000 & 1.4662 \\ 
321 & 23 & 225.5477 & 224.2384 & 2.4782  & 
322 & 23 & 226.5406 & 225.0000 & 2.9091  & 
323 & 23 & 227.5362 & 226.0000 & 2.9004 \\ 
324 & 23 & 228.5292 & 227.0000 & 2.8862  & 
325 & 23 & 229.5248 & 228.0000 & 2.8775  & 
326 & 23 & 230.5178 & 229.0000 & 2.8635 \\ 
327 & 23 & 231.5134 & 230.0000 & 2.8549  & 
328 & 23 & 232.5064 & 231.0000 & 2.8411  & 
329 & 23 & 233.5021 & 231.2384 & 4.8021 \\ 
 \hline 330 & 23 & 234.4951 & 232.0000 & 5.6377  & 
331 & 23 & 235.4908 & 233.0000 & 5.6208  & 
332 & 23 & 236.2246 & 234.0000 & 4.6738 \\ 
333 & 23 & 237.2203 & 235.0000 & 4.6600  & 
334 & 23 & 238.2135 & 236.0000 & 4.6381  & 
335 & 23 & 239.2093 & 237.0000 & 4.6244 \\ 
336 & 23 & 240.2025 & 238.0000 & 4.6027  & 
337 & 23 & 241.1982 & 239.0000 & 4.5892  & 
338 & 23 & 241.9368 & 240.0000 & 3.8287 \\ 
339 & 23 & 242.9326 & 241.0000 & 3.8175  & 
340 & 23 & 243.9260 & 242.0000 & 3.8001  & 
341 & 23 & 244.9218 & 243.0000 & 3.7890 \\ 
342 & 23 & 245.9153 & 244.0000 & 3.7718  & 
343 & 23 & 246.9111 & 245.0000 & 3.7609  & 
344 & 23 & 247.7373 & 246.0000 & 3.3340 \\ 
 \hline 345 & 23 & 248.7331 & 247.0000 & 3.3244  & 
346 & 23 & 249.7266 & 248.0000 & 3.3096  & 
347 & 23 & 250.7225 & 249.0000 & 3.3001 \\ 
348 & 23 & 251.6332 & 250.0000 & 3.1020  & 
349 & 23 & 252.6291 & 251.0000 & 3.0932  & 
350 & 23 & 253.4583 & 252.0000 & 2.7479 \\ 
351 & 23 & 254.4542 & 253.0000 & 2.7401  & 
352 & 23 & 255.4479 & 254.0000 & 2.7282  & 
353 & 23 & 256.4439 & 255.0000 & 2.7205 \\ 
354 & 23 & 257.1944 & 256.0000 & 2.2885  & 
355 & 23 & 258.1904 & 257.0000 & 2.2822  & 
356 & 23 & 259.1843 & 258.0000 & 2.2725 \\ 
357 & 23 & 260.1803 & 259.0000 & 2.2662  & 
358 & 23 & 261.1741 & 260.0000 & 2.2566  & 
359 & 23 & 262.1701 & 261.0000 & 2.2503 \\ 
 \hline 360 & 23 & 263.0839 & 262.0000 & 2.1198  & 
361 & 23 & 264.0799 & 263.0000 & 2.1140  & 
362 & 23 & 264.8360 & 264.0000 & 1.7852 \\ 
363 & 23 & 265.8321 & 265.0000 & 1.7803  & 
364 & 23 & 266.8262 & 266.0000 & 1.7730  & 
365 & 23 & 267.8223 & 267.0000 & 1.7682 \\ 
366 & 23 & 268.8164 & 268.0000 & 1.7610  & 
367 & 23 & 269.8125 & 269.0000 & 1.7562  & 
368 & 23 & 270.5726 & 270.0000 & 1.4872 \\ 
369 & 23 & 271.5687 & 271.0000 & 1.4832  & 
370 & 23 & 272.5630 & 272.0000 & 1.4773  & 
371 & 23 & 273.5591 & 273.0000 & 1.4733 \\ 
372 & 23 & 274.5534 & 274.0000 & 1.4675  & 
373 & 23 & 275.5495 & 275.0000 & 1.4636  & 
374 & 23 & 276.3135 & 276.0000 & 1.2427 \\ 
 \hline 375 & 23 & 277.3097 & 277.0000 & 1.2395  & 
376 & 23 & 278.3041 & 278.0000 & 1.2346  & 
377 & 23 & 279.3003 & 279.0000 & 1.2314 \\ 
378 & 23 & 280.2947 & 280.0000 & 1.2266  & 
379 & 23 & 281.2909 & 281.0000 & 1.2234  & 
380 & 23 & 282.1339 & 282.0000 & 1.0972 \\ 
381 & 23 & 283.1301 & 283.0000 & 1.0944  & 
382 & 23 & 284.1246 & 284.0000 & 1.0902  & 
383 & 23 & 285.1209 & 285.0000 & 1.0874 \\ 
384 & 23 & 285.8911 & 285.0000 & 1.8546  & 
385 & 23 & 286.8874 & 286.0000 & 1.8498  & 
386 & 23 & 287.8820 & 287.0000 & 1.8429 \\ 
387 & 23 & 288.8783 & 288.0000 & 1.8382  & 
388 & 23 & 289.8729 & 289.0000 & 1.8314  & 
389 & 23 & 290.8693 & 290.0000 & 1.8267 \\ 
 \hline 390 & 23 & 291.5701 & 291.0000 & 1.4846  & 
391 & 23 & 292.5665 & 292.0000 & 1.4809  & 
392 & 23 & 293.5612 & 293.0000 & 1.4755 \\ 
393 & 23 & 294.5576 & 294.0000 & 1.4718  & 
394 & 23 & 295.5524 & 295.0000 & 1.4665  & 
395 & 23 & 296.5488 & 296.0000 & 1.4628 \\ 
396 & 23 & 297.5436 & 297.0000 & 1.4576  & 
397 & 23 & 298.5400 & 298.0000 & 1.4539  & 
398 & 23 & 299.3902 & 299.0000 & 1.3106 \\ 
399 & 23 & 300.3866 & 300.0000 & 1.3073  & 
400 & 23 & 301.3815 & 301.0000 & 1.3027  & 
401 & 23 & 302.3779 & 302.0000 & 1.2995 \\ 
402 & 23 & 303.0877 & 303.0000 & 1.0627  & 
403 & 23 & 304.0842 & 304.0000 & 1.0601  & 
404 & 23 & 305.0792 & 305.0000 & 1.0565 \\ 
 \hline 405 & 23 & 306.0757 & 306.0000 & 1.0539  & 
406 & 23 & 307.0708 & 307.0000 & 1.0503  & 
407 & 23 & 308.0672 & 308.0000 & 1.0477 \\ 
408 & 23 & 309.0623 & 309.0000 & 1.0441  & 
409 & 23 & 310.0588 & 310.0000 & 1.0416  & 
305 & 25 & 202.4382 & 202.2384 & 1.1485 \\ 
306 & 25 & 203.4298 & 203.0000 & 1.3470  & 
307 & 25 & 204.4250 & 204.0000 & 1.3426  & 
308 & 25 & 205.2113 & 205.0000 & 1.1577 \\ 
309 & 25 & 206.2066 & 206.0000 & 1.1540  & 
310 & 25 & 207.1984 & 207.0000 & 1.1474  & 
311 & 25 & 208.1937 & 208.0000 & 1.1437 \\ 
312 & 25 & 209.0838 & 209.0000 & 1.0598  & 
313 & 25 & 210.0792 & 209.2384 & 1.7910  & 
314 & 25 & 210.8696 & 210.0000 & 1.8271 \\ 
 \hline \hline 
\end{tabular}

\begin{tabular}{||c|c|c|c|c||c|c|c|c|c||c|c|c|c|c||}
\hline\hline
 n & d & new & old & ratio &
n & d & new & old & ratio & n & d & new & old & ratio \\
\hline\hline
315 & 25 & 211.8650 & 211.0000 & 1.8214  & 
316 & 25 & 212.8570 & 212.0000 & 1.8113  & 
317 & 25 & 213.8524 & 213.0000 & 1.8055 \\ 
318 & 25 & 214.1577 & 214.0000 & 1.1155  & 
319 & 25 & 215.1533 & 215.0000 & 1.1121  & 
320 & 25 & 216.1456 & 216.0000 & 1.1062 \\ 
321 & 25 & 217.1412 & 216.2384 & 1.8697  & 
322 & 25 & 218.1335 & 217.0000 & 2.1939  & 
323 & 25 & 219.1291 & 218.0000 & 2.1872 \\ 
324 & 25 & 220.1215 & 219.0000 & 2.1757  & 
325 & 25 & 221.1171 & 220.0000 & 2.1691  & 
326 & 25 & 222.1095 & 221.0000 & 2.1577 \\ 
327 & 25 & 223.1051 & 222.0000 & 2.1511  & 
328 & 25 & 224.0975 & 223.0000 & 2.1399  & 
329 & 25 & 225.0932 & 223.2384 & 3.6169 \\ 
 \hline 330 & 25 & 226.0856 & 224.0000 & 4.2445  & 
331 & 25 & 227.0812 & 225.0000 & 4.2317  & 
332 & 25 & 227.7886 & 226.0000 & 3.4547 \\ 
333 & 25 & 228.7843 & 227.0000 & 3.4445  & 
334 & 25 & 229.7769 & 228.0000 & 3.4269  & 
335 & 25 & 230.7726 & 229.0000 & 3.4167 \\ 
336 & 25 & 231.7653 & 230.0000 & 3.3994  & 
337 & 25 & 232.7610 & 230.2384 & 5.7460  & 
338 & 25 & 233.4735 & 231.0000 & 5.5540 \\ 
339 & 25 & 234.4693 & 232.0000 & 5.5378  & 
340 & 25 & 235.4621 & 233.0000 & 5.5103  & 
341 & 25 & 236.4579 & 234.0000 & 5.4942 \\ 
342 & 25 & 237.4508 & 235.0000 & 5.4670  & 
343 & 25 & 238.4465 & 236.0000 & 5.4510  & 
344 & 25 & 239.2554 & 237.0000 & 4.7747 \\ 
 \hline 345 & 25 & 240.2513 & 238.0000 & 4.7610  & 
346 & 25 & 241.2442 & 239.0000 & 4.7378  & 
347 & 25 & 242.2400 & 240.0000 & 4.7241 \\ 
348 & 25 & 243.1419 & 241.0000 & 4.4134  & 
349 & 25 & 244.1377 & 242.0000 & 4.4007  & 
350 & 25 & 244.9499 & 243.0000 & 3.8636 \\ 
351 & 25 & 245.9459 & 244.0000 & 3.8527  & 
352 & 25 & 246.9390 & 245.0000 & 3.8344  & 
353 & 25 & 247.9349 & 246.0000 & 3.8235 \\ 
354 & 25 & 248.6605 & 247.0000 & 3.1614  & 
355 & 25 & 249.6565 & 248.0000 & 3.1525  & 
356 & 25 & 250.6498 & 249.0000 & 3.1379 \\ 
357 & 25 & 251.6458 & 250.0000 & 3.1292  & 
358 & 25 & 252.6391 & 251.0000 & 3.1147  & 
359 & 25 & 253.6351 & 252.0000 & 3.1060 \\ 
 \hline 360 & 25 & 254.5403 & 253.0000 & 2.9085  & 
361 & 25 & 255.5363 & 254.0000 & 2.9004  & 
362 & 25 & 256.2680 & 255.0000 & 2.4083 \\ 
363 & 25 & 257.2641 & 256.0000 & 2.4017  & 
364 & 25 & 258.2576 & 257.0000 & 2.3910  & 
365 & 25 & 259.2537 & 258.0000 & 2.3844 \\ 
366 & 25 & 260.2472 & 259.0000 & 2.3738  & 
367 & 25 & 261.2432 & 260.0000 & 2.3673  & 
368 & 25 & 261.9794 & 261.0000 & 1.9717 \\ 
369 & 25 & 262.9755 & 262.0000 & 1.9664  & 
370 & 25 & 263.9692 & 263.0000 & 1.9578  & 
371 & 25 & 264.9653 & 264.0000 & 1.9525 \\ 
372 & 25 & 265.9590 & 265.0000 & 1.9440  & 
373 & 25 & 266.9551 & 266.0000 & 1.9388  & 
374 & 25 & 267.6956 & 267.0000 & 1.6195 \\ 
 \hline 375 & 25 & 268.6918 & 268.0000 & 1.6153  & 
376 & 25 & 269.6856 & 269.0000 & 1.6084  & 
377 & 25 & 270.6818 & 270.0000 & 1.6041 \\ 
378 & 25 & 271.6756 & 271.0000 & 1.5973  & 
379 & 25 & 272.6718 & 272.0000 & 1.5931  & 
380 & 25 & 273.4991 & 273.0000 & 1.4133 \\ 
381 & 25 & 274.4953 & 274.0000 & 1.4096  & 
382 & 25 & 275.4893 & 275.0000 & 1.4037  & 
383 & 25 & 276.4855 & 276.0000 & 1.4001 \\ 
384 & 25 & 277.2327 & 276.0000 & 2.3501  & 
385 & 25 & 278.2290 & 277.0000 & 2.3441  & 
386 & 25 & 279.2231 & 278.0000 & 2.3345 \\ 
387 & 25 & 280.2194 & 279.0000 & 2.3285  & 
388 & 25 & 281.2135 & 280.0000 & 2.3189  & 
389 & 25 & 282.2097 & 281.0000 & 2.3130 \\ 
 \hline 390 & 25 & 282.8807 & 282.0000 & 1.8412  & 
391 & 25 & 283.8770 & 283.0000 & 1.8366  & 
392 & 25 & 284.8713 & 284.0000 & 1.8293 \\ 
393 & 25 & 285.8676 & 285.0000 & 1.8247  & 
394 & 25 & 286.8619 & 286.0000 & 1.8174  & 
395 & 25 & 287.8583 & 287.0000 & 1.8128 \\ 
396 & 25 & 288.8525 & 288.0000 & 1.8057  & 
397 & 25 & 289.8489 & 289.0000 & 1.8011  & 
398 & 25 & 290.6841 & 290.0000 & 1.6067 \\ 
399 & 25 & 291.6805 & 291.0000 & 1.6027  & 
400 & 25 & 292.6749 & 292.0000 & 1.5965  & 
401 & 25 & 293.6713 & 293.0000 & 1.5925 \\ 
402 & 25 & 294.3521 & 294.0000 & 1.2764  & 
403 & 25 & 295.3485 & 295.0000 & 1.2733  & 
404 & 25 & 296.3430 & 296.0000 & 1.2684 \\ 
 \hline 405 & 25 & 297.3395 & 297.0000 & 1.2653  & 
406 & 25 & 298.3340 & 298.0000 & 1.2605  & 
407 & 25 & 299.3305 & 299.0000 & 1.2574 \\ 
408 & 25 & 300.3250 & 300.0000 & 1.2527  & 
409 & 25 & 301.3215 & 301.0000 & 1.2496  & 
306 & 27 & 195.1500 & 195.0000 & 1.1095 \\ 
307 & 27 & 196.1453 & 196.0000 & 1.1060  & 
313 & 27 & 201.7703 & 201.2384 & 1.4458  & 
314 & 27 & 202.5418 & 202.0000 & 1.4558 \\ 
315 & 27 & 203.5373 & 203.0000 & 1.4513  & 
316 & 27 & 204.5288 & 204.0000 & 1.4427  & 
317 & 27 & 205.5243 & 205.0000 & 1.4382 \\ 
321 & 27 & 208.7498 & 208.2384 & 1.4255  & 
322 & 27 & 209.7417 & 209.0000 & 1.6721  & 
323 & 27 & 210.7373 & 210.0000 & 1.6671 \\ 
 \hline 324 & 27 & 211.7292 & 211.0000 & 1.6577  & 
325 & 27 & 212.7248 & 212.0000 & 1.6527  & 
326 & 27 & 213.7168 & 213.0000 & 1.6435 \\ 
327 & 27 & 214.7124 & 214.0000 & 1.6385  & 
328 & 27 & 215.7044 & 215.0000 & 1.6294  & 
329 & 27 & 216.7000 & 215.2384 & 2.7541 \\ 
330 & 27 & 217.6920 & 216.0000 & 3.2310  & 
331 & 27 & 218.6877 & 217.0000 & 3.2213  & 
332 & 27 & 219.3686 & 218.0000 & 2.5822 \\ 
333 & 27 & 220.3643 & 219.0000 & 2.5746  & 
334 & 27 & 221.3565 & 220.0000 & 2.5606  & 
335 & 27 & 222.3522 & 221.0000 & 2.5530 \\ 
336 & 27 & 223.3444 & 222.0000 & 2.5392  & 
337 & 27 & 224.3401 & 222.2384 & 4.2921  & 
338 & 27 & 225.0267 & 223.0000 & 4.0748 \\ 
 \hline 339 & 27 & 226.0225 & 224.0000 & 4.0630  & 
340 & 27 & 227.0149 & 225.0000 & 4.0414  & 
341 & 27 & 228.0107 & 226.0000 & 4.0297 \\ 
342 & 27 & 229.0030 & 227.0000 & 4.0084  & 
343 & 27 & 229.9988 & 228.0000 & 3.9968  & 
344 & 27 & 230.7905 & 229.0000 & 3.4594 \\ 
345 & 27 & 231.7864 & 229.2384 & 5.8480  & 
346 & 27 & 232.7788 & 230.0000 & 6.8630  & 
347 & 27 & 233.7747 & 231.0000 & 6.8433 \\ 
348 & 27 & 234.6677 & 232.0000 & 6.3543  & 
349 & 27 & 235.6636 & 233.0000 & 6.3362  & 
350 & 27 & 236.4589 & 234.0000 & 5.4980 \\ 
351 & 27 & 237.4548 & 235.0000 & 5.4824  & 
352 & 27 & 238.4475 & 236.0000 & 5.4546  & 
353 & 27 & 239.4434 & 237.0000 & 5.4392 \\ 
 \hline 354 & 27 & 240.1442 & 238.0000 & 4.4206  & 
355 & 27 & 241.1402 & 239.0000 & 4.4083  & 
356 & 27 & 242.1330 & 240.0000 & 4.3863 \\ 
357 & 27 & 243.1290 & 241.0000 & 4.3742  & 
358 & 27 & 244.1218 & 242.0000 & 4.3525  & 
359 & 27 & 245.1178 & 243.0000 & 4.3404 \\ 
360 & 27 & 246.0145 & 244.0000 & 4.0405  & 
361 & 27 & 247.0105 & 245.0000 & 4.0293  & 
362 & 27 & 247.7180 & 246.0000 & 3.2898 \\ 
363 & 27 & 248.7141 & 247.0000 & 3.2809  & 
364 & 27 & 249.7071 & 248.0000 & 3.2651  & 
365 & 27 & 250.7032 & 249.0000 & 3.2562 \\ 
366 & 27 & 251.6962 & 250.0000 & 3.2405  & 
367 & 27 & 252.6923 & 251.0000 & 3.2317  & 
368 & 27 & 253.4046 & 252.0000 & 2.6474 \\ 
 \hline 369 & 27 & 254.4007 & 253.0000 & 2.6403  & 
370 & 27 & 255.3939 & 254.0000 & 2.6279  & 
371 & 27 & 256.3900 & 255.0000 & 2.6209 \\ 
372 & 27 & 257.3832 & 256.0000 & 2.6085  & 
373 & 27 & 258.3794 & 257.0000 & 2.6016  & 
374 & 27 & 259.0963 & 258.0000 & 2.1381 \\ 
375 & 27 & 260.0925 & 259.0000 & 2.1324  & 
376 & 27 & 261.0858 & 260.0000 & 2.1226  & 
377 & 27 & 262.0820 & 261.0000 & 2.1170 \\ 
378 & 27 & 263.0754 & 262.0000 & 2.1073  & 
379 & 27 & 264.0716 & 263.0000 & 2.1017  & 
380 & 27 & 264.8832 & 264.0000 & 1.8445 \\ 
381 & 27 & 265.8795 & 265.0000 & 1.8397  & 
382 & 27 & 266.8729 & 266.0000 & 1.8314  & 
383 & 27 & 267.8692 & 267.0000 & 1.8266 \\ 
 \hline 384 & 27 & 268.5935 & 267.0000 & 3.0178  & 
385 & 27 & 269.5898 & 268.0000 & 3.0100  & 
386 & 27 & 270.5834 & 269.0000 & 2.9967 \\ 
387 & 27 & 271.5797 & 270.0000 & 2.9890  & 
388 & 27 & 272.5733 & 271.0000 & 2.9758  & 
389 & 27 & 273.5696 & 272.0000 & 2.9682 \\ 
390 & 27 & 274.2106 & 273.0000 & 2.3144  & 
391 & 27 & 275.2070 & 274.0000 & 2.3086  & 
392 & 27 & 276.2007 & 275.0000 & 2.2986 \\ 
393 & 27 & 277.1971 & 276.0000 & 2.2928  & 
394 & 27 & 278.1909 & 277.0000 & 2.2829  & 
395 & 27 & 279.1873 & 278.0000 & 2.2772 \\ 
396 & 27 & 280.1810 & 279.0000 & 2.2674  & 
397 & 27 & 281.1774 & 280.0000 & 2.2617  & 
398 & 27 & 281.9977 & 281.0000 & 1.9968 \\ 
 \hline \hline 
\end{tabular}

\begin{tabular}{||c|c|c|c|c||c|c|c|c|c||c|c|c|c|c||}
\hline\hline
 n & d & new & old & ratio &
n & d & new & old & ratio & n & d & new & old & ratio \\
\hline\hline
399 & 27 & 282.9941 & 282.0000 & 1.9918  & 
400 & 27 & 283.9880 & 283.0000 & 1.9834  & 
401 & 27 & 284.9844 & 284.0000 & 1.9785 \\ 
402 & 27 & 285.6362 & 285.0000 & 1.5542  & 
403 & 27 & 286.6327 & 286.0000 & 1.5504  & 
404 & 27 & 287.6267 & 287.0000 & 1.5440 \\ 
405 & 27 & 288.6232 & 288.0000 & 1.5402  & 
406 & 27 & 289.6172 & 289.0000 & 1.5339  & 
407 & 27 & 290.6137 & 290.0000 & 1.5302 \\ 
408 & 27 & 291.6077 & 291.0000 & 1.5239  & 
409 & 27 & 292.6042 & 292.0000 & 1.5202  & 
410 & 27 & 293.1799 & 293.0000 & 1.1328 \\ 
411 & 27 & 294.1764 & 294.0000 & 1.1301  & 
412 & 27 & 295.1706 & 295.0000 & 1.1256  & 
413 & 27 & 296.1672 & 296.0000 & 1.1229 \\ 
 \hline 414 & 27 & 297.1614 & 297.0000 & 1.1184  & 
415 & 27 & 298.1580 & 298.0000 & 1.1157  & 
416 & 27 & 299.1522 & 299.0000 & 1.1113 \\ 
417 & 27 & 300.1488 & 300.0000 & 1.1087  & 
418 & 27 & 301.1431 & 301.0000 & 1.1043  & 
419 & 27 & 302.1397 & 302.0000 & 1.1016 \\ 
420 & 27 & 303.0515 & 303.0000 & 1.0364  & 
421 & 27 & 304.0481 & 304.0000 & 1.0339  & 
314 & 29 & 194.2232 & 194.0000 & 1.1673 \\ 
315 & 29 & 195.2186 & 195.0000 & 1.1636  & 
316 & 29 & 196.2099 & 196.0000 & 1.1566  & 
317 & 29 & 197.2053 & 197.0000 & 1.1529 \\ 
321 & 29 & 200.3679 & 200.2384 & 1.0939  & 
322 & 29 & 201.3595 & 201.0000 & 1.2830  & 
323 & 29 & 202.3551 & 202.0000 & 1.2791 \\ 
 \hline 324 & 29 & 203.3467 & 203.0000 & 1.2717  & 
325 & 29 & 204.3424 & 204.0000 & 1.2678  & 
326 & 29 & 205.3340 & 205.0000 & 1.2605 \\ 
327 & 29 & 206.3296 & 206.0000 & 1.2567  & 
328 & 29 & 207.3213 & 207.0000 & 1.2495  & 
329 & 29 & 208.3169 & 207.2384 & 2.1119 \\ 
330 & 29 & 209.3086 & 208.0000 & 2.4771  & 
331 & 29 & 210.3043 & 209.0000 & 2.4696  & 
332 & 29 & 210.9590 & 210.0000 & 1.9440 \\ 
333 & 29 & 211.9547 & 211.0000 & 1.9382  & 
334 & 29 & 212.9466 & 212.0000 & 1.9273  & 
335 & 29 & 213.9423 & 213.0000 & 1.9216 \\ 
336 & 29 & 214.9342 & 214.0000 & 1.9108  & 
337 & 29 & 215.9299 & 214.2384 & 3.2298  & 
338 & 29 & 216.5907 & 215.0000 & 3.0120 \\ 
 \hline 339 & 29 & 217.5865 & 216.0000 & 3.0033  & 
340 & 29 & 218.5786 & 217.0000 & 2.9867  & 
341 & 29 & 219.5743 & 218.0000 & 2.9780 \\ 
342 & 29 & 220.5664 & 219.0000 & 2.9617  & 
343 & 29 & 221.5622 & 220.0000 & 2.9530  & 
344 & 29 & 222.3368 & 221.0000 & 2.5260 \\ 
345 & 29 & 223.3327 & 221.2384 & 4.2701  & 
346 & 29 & 224.3248 & 222.0000 & 5.0101  & 
347 & 29 & 225.3207 & 223.0000 & 4.9956 \\ 
348 & 29 & 226.2051 & 224.0000 & 4.6111  & 
349 & 29 & 227.2010 & 225.0000 & 4.5979  & 
350 & 29 & 227.9795 & 226.0000 & 3.9436 \\ 
351 & 29 & 228.9754 & 227.0000 & 3.9324  & 
352 & 29 & 229.9677 & 228.0000 & 3.9116  & 
353 & 29 & 230.9636 & 228.2384 & 6.6127 \\ 
 \hline 354 & 29 & 231.6399 & 229.0000 & 6.2327  & 
355 & 29 & 232.6358 & 230.0000 & 6.2153  & 
356 & 29 & 233.6283 & 231.0000 & 6.1830 \\ 
357 & 29 & 234.6243 & 232.0000 & 6.1657  & 
358 & 29 & 235.6168 & 233.0000 & 6.1338  & 
359 & 29 & 236.6127 & 234.0000 & 6.1166 \\ 
360 & 29 & 237.5011 & 235.0000 & 5.6612  & 
361 & 29 & 238.4971 & 236.0000 & 5.6455  & 
362 & 29 & 239.1805 & 237.0000 & 4.5330 \\ 
363 & 29 & 240.1765 & 238.0000 & 4.5207  & 
364 & 29 & 241.1692 & 239.0000 & 4.4978  & 
365 & 29 & 242.1653 & 240.0000 & 4.4855 \\ 
366 & 29 & 243.1580 & 241.0000 & 4.4629  & 
367 & 29 & 244.1540 & 242.0000 & 4.4507  & 
368 & 29 & 244.8426 & 243.0000 & 3.5865 \\ 
 \hline 369 & 29 & 245.8387 & 244.0000 & 3.5769  & 
370 & 29 & 246.8316 & 245.0000 & 3.5592  & 
371 & 29 & 247.8277 & 246.0000 & 3.5497 \\ 
372 & 29 & 248.8205 & 247.0000 & 3.5321  & 
373 & 29 & 249.8167 & 248.0000 & 3.5226  & 
374 & 29 & 250.5102 & 249.0000 & 2.8485 \\ 
375 & 29 & 251.5064 & 250.0000 & 2.8410  & 
376 & 29 & 252.4994 & 251.0000 & 2.8272  & 
377 & 29 & 253.4956 & 252.0000 & 2.8198 \\ 
378 & 29 & 254.4886 & 253.0000 & 2.8061  & 
379 & 29 & 255.4848 & 254.0000 & 2.7987  & 
380 & 29 & 256.2809 & 255.0000 & 2.4299 \\ 
381 & 29 & 257.2771 & 256.0000 & 2.4236  & 
382 & 29 & 258.2702 & 257.0000 & 2.4120  & 
383 & 29 & 259.2665 & 258.0000 & 2.4057 \\ 
 \hline 384 & 29 & 259.9680 & 258.0000 & 3.9123  & 
385 & 29 & 260.9643 & 259.0000 & 3.9022  & 
386 & 29 & 261.9575 & 260.0000 & 3.8840 \\ 
387 & 29 & 262.9538 & 261.0000 & 3.8740  & 
388 & 29 & 263.9471 & 262.0000 & 3.8559  & 
389 & 29 & 264.9433 & 263.0000 & 3.8459 \\ 
390 & 29 & 265.5547 & 264.0000 & 2.9377  & 
391 & 29 & 266.5510 & 265.0000 & 2.9303  & 
392 & 29 & 267.5444 & 266.0000 & 2.9169 \\ 
393 & 29 & 268.5408 & 267.0000 & 2.9095  & 
394 & 29 & 269.5342 & 268.0000 & 2.8963  & 
395 & 29 & 270.5306 & 269.0000 & 2.8890 \\ 
396 & 29 & 271.5240 & 270.0000 & 2.8759  & 
397 & 29 & 272.5203 & 271.0000 & 2.8686  & 
398 & 29 & 273.3258 & 272.0000 & 2.5067 \\ 
 \hline 399 & 29 & 274.3222 & 273.0000 & 2.5004  & 
400 & 29 & 275.3157 & 274.0000 & 2.4892  & 
401 & 29 & 276.3121 & 275.0000 & 2.4830 \\ 
402 & 29 & 276.9350 & 276.0000 & 1.9119  & 
403 & 29 & 277.9315 & 277.0000 & 1.9073  & 
404 & 29 & 278.9251 & 278.0000 & 1.8989 \\ 
405 & 29 & 279.9216 & 279.0000 & 1.8942  & 
406 & 29 & 280.9153 & 280.0000 & 1.8859  & 
407 & 29 & 281.9118 & 281.0000 & 1.8813 \\ 
408 & 29 & 282.9055 & 282.0000 & 1.8731  & 
409 & 29 & 283.9019 & 283.0000 & 1.8685  & 
410 & 29 & 284.4424 & 284.0000 & 1.3588 \\ 
411 & 29 & 285.4389 & 285.0000 & 1.3556  & 
412 & 29 & 286.4327 & 286.0000 & 1.3498  & 
413 & 29 & 287.4293 & 287.0000 & 1.3466 \\ 
 \hline 414 & 29 & 288.4231 & 288.0000 & 1.3409  & 
415 & 29 & 289.4197 & 289.0000 & 1.3377  & 
416 & 29 & 290.4136 & 290.0000 & 1.3320 \\ 
417 & 29 & 291.4101 & 291.0000 & 1.3288  & 
418 & 29 & 292.4040 & 292.0000 & 1.3232  & 
419 & 29 & 293.4006 & 293.0000 & 1.3200 \\ 
420 & 29 & 294.3052 & 294.0000 & 1.2356  & 
421 & 29 & 295.3018 & 295.0000 & 1.2327  & 
& & & & \\
\hline\hline
\end{tabular}

}


\begin{thebibliography}{20}

\bibitem{BC}
R.~C.~Bose and S.~Chowla, 
{\it Theorems in the additive theory of numbers},
Comment. Math. Helv. {\bf 37} (1962/1963), 141--147. 




\bibitem{Gil} E.~N.~Gilbert, {\it A comparison of signalling alphabets},
Bell Syst.~Tech.~J.~{\bf 31} (1952).







\bibitem{GS} R.~L.~Graham and N.~J.~A.~Sloane, {\it 
Lower bounds for constant weight codes}, IEEE Trans. Inform. Theory {\bf 26} 
(1980), no. 1, 37--43.

\bibitem{HR} H.~Halberstam and K.~F.~Roth, {\it Sequences}, Vol.~1, 
Oxford University Press, 1966.

\bibitem{Ham} R.~W.~Hamming, {\it Error detecting and error correcting codes},
Bell Syst.\ Tech.\ J.~{\bf 29} (1950), 147--160.





\bibitem{Li} S.~Litsyn, {\it An updated table of the best binary codes known}, 
Handbook of coding theory, Vol. I, II, 463--498, North-Holland, Amsterdam,
1998.

\bibitem{LRS}S.~Litsyn, E.~M.~Rains and 
N.~J.~A.~Sloane,\\ 
{\tt http://www.eng.tau.ac.il/\~{}litsyn/tableand/index.html}, or\\
{\tt http://www.research.att.com/\~{}njas/codes/And/}.

\bibitem{McWS}
F.~J.~MacWilliams and N.~J.~A.~Sloane, {\it 
The Theory of Error-Correcting Codes}, North-Holland Mathematical Library, 
Vol. {\bf 16}. North-Holland
Publishing Co., Amsterdam-New York-Oxford, 1977.

\bibitem{RS}E.~M.~Rains and 
N.~J.~A.~Sloane,\\ 
{\tt http://www.research.att.com/\~{}njas/codes/Andw/}.


\bibitem{Sacks}
G.~E.~Sacks, {\it Multiple error correction by means of parity checks},
IEEE Trans.\ Info.\ Theory {\bf 4} (1958), 145--147.

\bibitem{Var}
R.~R.~Varshamov, {\it Estimate of the number of signals in error correcting
codes}, Dokl.\ Akad.\ Nauk USSR {\bf 117} (1957), 739--741.


\end{thebibliography}
\end{document}